\documentclass[twoside, 12pt,amstex]{article}
\usepackage{amsmath,amssymb, soul,epsfig}
\usepackage{amsfonts}
\usepackage{amsthm}
\usepackage{mathrsfs}
\usepackage{enumerate}
\usepackage{amsthm}
\usepackage{graphicx}
\usepackage[english]{babel}
\usepackage[cp1251]{inputenc}
\date{\empty}

\numberwithin{equation}{section}

\def\@evenhead{\vbox{\hbox to \textwidth{\thepage\hfil\sl\leftmark\strut}\hrule}}
\def\@oddhead{\vbox{\hbox to \textwidth{\rightmark\hfill\thepage\strut}\hrule}}
\begin{document}
\sloppy

\centerline{\bf NEW 2-MICROLOCAL BESOV AND TRIEBEL--LIZORKIN SPACES} 
\centerline{\bf VIA THE LITLLEWOOD -- PALEY DECOMPOSITION}

\vskip 0.3cm

\centerline{\bf K. Saka}
\markboth{\hfill{\footnotesize\rm Koichi~Saka}\hfill}
{\hfill{\footnotesize\sl  New 2-microlocal Besov and Triebel-Lizorkin spaces via the Littlewood - Paley decomposition}\hfill}
\vskip0.3cm

\vskip0.7cm

\noindent
{\bf Key words:} 
wavelet,  Besov space,  Triebel--Lizorkin space,  pseudo-differential operator,  Calder$\acute{\rm o}$n--Zygmund operator,  atomic and molecular 
decomposition, 2-microlocal space,  $\varphi$--transform.

\vskip 0.2cm

\noindent
{\bf AMS Mathematics Subject Classification:} 
42B35, 42B20, 42B25, 42C40.

\vskip 0.2cm

\noindent
{\bf Abstract}. \ \ 
In this paper we introduce and investigate new 2-microlocal Besov and Triebel--Lizorkin spaces via the Littlewood - Paley decomposition.  
 We establish characterizations of these  function spaces by the $\varphi$--transform, the atomic and molecular decomposition  and the wavelet decomposition. 
 As applications we prove  boundedness of the 
the Calder$\acute{{\rm o}}$n--Zygmund operators  and the pseudo--differential operators on  the function spaces. Moreover, we give characterizations via oscillations and differences.

\section{\large Introduction}

It is well known that function spaces have increasing applications in  
many areas of modern analysis, in particular,  harmonic analysis and 
partial differential equations. 
The most general function spaces, probably, 
 are the  Besov spaces and  the Triebel--Lizorkin spaces which cover many 
 classical concrete function spaces such as Lebesgue spaces, 
 Lipschitz spaces, Sobolev spaces, Hardy spaces and  BMO spaces ([37], [38]).

 D. Yang and W. Yuan in [41], [42] and W. Sickel, D. Yang and W. Yuan in [36], 
  introduced a class of Besov type and Triebel--Lizorkin type spaces 
  which generalized many classical function spaces such as Besov spaces, Triebel-Lizorkin spaces, Morrey spaces and $Q$-type spaces. Recently  the Besov type and Triebel-Lizorkin type space with variable exponents was investigated by many authors (e.g. [43], [44]).   

The 2-microloal space is due to Bony [3] in order to study the propagation of singularities of the solutions of nonlinear evolution equations. 
It is an appropriate instrument to describe the local regularity and the oscillatory behavior of functions near  singularity (Meyer [32]). The theory has been elaborated and widely used in fractal analysis and signal processing.
 For systematic discussions of the concept and further references of 2-microlocal spaces, we refer to Meyer[31], [32],  Levy-Vehel and Seuret [30], 
Jaffard ([17], [18], [19], [20]), Jaffard and M${\rm \acute{e}}$lot [21], and  Jaffard and Meyer [22]. 
 
 The 2-microlocal spaces have been generalized by  Jaffard  as a general pointwise regularity associated with Banach  or quasi-Banach spaces [19], [20]. 
In this paper we introduce new inhomogeneous 2-microlocal spaces based 
on Jaffard's idea (See [33] for the homogeneous 2-microlocal spaces) and   
 we will investigate the properties and the characterizations of these 
new 2-microlocal  Besov  and Triebel--Lizorkin  spaces which unify many classical function spaces such as the Besov type and Triebel--Lizorkin type spaces, the 2--microlocal spaces in the sense of Meyer [32],  the Morrey space and  the local Morrey spaces. These new function space are very similar to the classical 2-microlocal Besov and Triebel-Lizorkin spaces  studied recently by many authors ([1], [6], [8], [13], [14], [15], [16], [25], [26], [39], [40]). 

The plan of the remaining sections in the paper is as follows:

In Section 2 we give the definitions of our  new 2-microlocal  spaces  via the 
Littlewood-Paley decomposition  and the 
notations which are used later and  we give examples for these spaces.
 
In Section 3 we define corresponding sequence spaces 
for our function spaces. 
Furthermore, 
we give  some auxiliary lemmas which are needed in
later sections.

In Section 4 we will characterize our function spaces 
via the corresponding sequence spaces
by the $\varphi$--transform in the sense of Fraizer--Jarwerth [10], the atomic and molecular decomposition and the wavelet decomposition.
Moreover, we investigate the properties for these function spaces and 
we also study relations between  
our 2-microlocal  spaces and the classical 2-microlocal spaces. 

 In Section 5, as  applications, we give  the conditions under 
which 
the Calder$\acute{{\rm o}}$n--Zygmund operators  and the pseudo--differential 
operators are bounded on the function spaces.

In Section 6 we give the characterizations via differences and oscillations. 

Throughout the paper, we use $C$ to denote a positive constant. 
But the same notation $C$ are 
not necessarily the same on any two occurrences.  
 We use the  notations $i\vee j=\max\{i,\ j\}$, $i\wedge j=\min\{i,\ j\}$, and\ $a_{+}=a\vee 0$.
The symbol $X \sim Y$ 
 means that there exist positive constants $C_{1}$ and $C_{2}$ such that 
$X \leq C_{1} Y$ and $Y \leq C_{2}X$.

\section{Definitions}
We consider the dyadic cubes in ${\Bbb R}^n$ of the form 
 $Q=[0,\  2^{-l})^n + 2^{-l}k$ for $k \in {\Bbb Z}^n$ and $ l \in {\Bbb Z}$, 
and use the notation $l(Q)=2^{-l}$ for the side length and $x_{Q}= 2^{-l}k$ 
for the corner point. Throughout the paper, we use the notations $P,\ Q,\ R\ $for the dyadic cubes of the form $[0,\  2^{-l})^n + 2^{-l}k$ in ${\mathbb R}^n$, and when the dyadic cubes $Q$ appear as indices, it is understood that $Q$ runs over all dyadic cubes of this form in ${\mathbb R}^n$. 
We denote by $\mathcal{D}$ the set of all dyadic cubes of this form. 
For a dyadic cube $Q$ and a constant $c>1$, $cQ$ denotes the cube of same center as  $Q$ and $c$ times larger.
We denote by $\chi_{E}$ the characteristic function of a set $E$ 
in ${\mathbb R}^n$. 

We set $\mathbb{N}= 
\{ 1, 2, \cdots \}$ and $\mathbb{N}_{0} = \mathbb{N} \cup \{ 0 \}$.
 Let ${\mathcal S}={\mathcal S}({\Bbb R}^n)$ be the space of 
all Schwartz functions on 
${\Bbb R}^n$ and ${\mathcal S}'$ its dual. 
 
We use $\langle f,\ g \rangle$ for the standard inner product $\int f \bar{g}$ 
of two functions and the same notation is employed for the action of a distribution $f \in \mathcal{S}'$
on $ \bar{g} \in \mathcal{S}$.

\noindent
Let $\phi_{0}$ be a Schwartz function and $\hat{\phi_{0}}$ its Fourier transform satisfying

(1.1)\ \  
supp $\hat{\phi_{0}} \subset 
\{ \xi \in {\Bbb R}^n: \  |\xi| \leq 2 \}$,

(1.2)\ \ $\hat{\phi_{0}}(\xi) = 1$ if $|\xi| \leq 1$.

We set

$\phi(x)=\phi_{0}(x)-2^{-n}\phi_{0}(2^{-1}x)$, 
$\phi_{0}^j=2^{jn}\phi_{0}(2^jx)$, $S_{j}f=f*\phi_{0}^j$  for $j \in \mathbb{N}_{0}$, and 
$\phi_{j}(x)=2^{jn}\phi(2^jx)$  for $j \in \mathbb{N}$.

Then we have 

(1.3)\ \ supp $\hat{\phi} \subset 
\{ \xi \in {\Bbb R}^n: \  \frac{1}{2} \leq |\xi| \leq 2 \}$, and 

(1.4)\ \ there exist positive numbers c and a sufficiently small $\epsilon$ such that  $\hat{\phi}(\xi) \geq c$ in $1-\epsilon \leq |\xi| \leq 1+\epsilon$. 

It holds that $\sum_{j \in \mathbb{N}_{0}}\hat{\phi}_{j}=1$.
Let $f \in {\mathcal S}'$, then  we have the Littlewood-Paley decomposition 
$f= \sum_{j \in \mathbb{N}_{0} }f*\phi_{j}$
 (convergence in ${\mathcal S}'$) [36, Triebel 2.3.1(6)].

 Let $s \in {\Bbb R}$. For $f \in {\mathcal S}'$, 
we define  some sequences indexed by dyadic cubes $P$:

\bigskip

$c(B^s_{pq})(P)
=(\sum_{i \geq (-\log_{2}l(P))\vee 0}||2^{is}
f*\phi_{i}||_{L^p(P)}^{q})^{1/q}$, 
 $0< p, q \leq \infty$,

$c(F^s_{pq})(P)=
||\{\sum_{i \geq (-\log_{2}l(P))\vee 0}(2^{is}
|f*\phi_{i}|)^q\}^{1/q}||_{L^p(P)}$, 

 $0< p < \infty$,\ $0 < q \leq \infty$,

$c(F^s_{\infty q})(P)=l(P)^{-\frac{n}{q}}
||\{\sum_{i \geq (-\log_{2}l(P))\vee 0}(2^{is}
|f*\phi_{i}|)^q\}^{1/q}||_{L^q(P)}$,

 $0 < q \leq \infty$,

\noindent
with the usual modification for $q=\infty$.

 We  shall use the notation $E^s_{pq}$ with either $B^s_{pq}$ or 
$F^s_{pq}$. 
We say the B-type case when $E^{s'}_{pq}=B^{s'}_{pq}$, and the F-type case when $E^{s'}_{pq}=F^{s'}_{pq}$.

\bigskip

\noindent
{\bf Definition 1.}\ \ 
 Let $s, \ s', \ \sigma\ \in {\Bbb R},\ 0 < p,q \leq \infty$ and 
 $x_{0} \in \mathbb{R}^n$.

The space $A^s(E^{s'}_{pq})^{\sigma}_{x_{0}}$ is defined to be 
the space of all $f \in {\mathcal S}'$ such that

$$
 ||f||_{A^s(E^{s'}_{pq})^{\sigma}_{x_{0}}}\equiv 
\sup_{ \mathcal{D} \ni Q \ni x_{0} }l(Q)^{-\sigma}
\sup_{ \mathcal{D} \ni P \subset 3Q\ }l(P)^{-s}
c(E^{s'}_{pq})(P) < \infty.
$$

\bigskip

\noindent
The following abbreviation 
$A^0(E^{s'}_{pq})^{\sigma}_{x_{0}} \equiv 
(E^{s'}_{pq})^{\sigma}_{x_{0}}$,
$A^s(E^{s'}_{pq})^{0}_{x_{0}} \equiv A^s(E^{s'}_{pq})$ and 
$A^0(E^{s'}_{pq})^{0}_{x_{0}}\equiv E^{s'}_{pq} 
\equiv E^{s'}_{pq}(\mathbb{R}^n)$ will be used in the sequel. 
We note that the space 
$A^s(E^{s'}_{pq})$ is the inhomogeneous Besov type space or 
the inhomogeneous Triebel--Lizorkin type space  in the sense of Yang--Sickel--Yuan [26] and  the space $E^{s'}_{pq}\equiv E^{s'}_{pq}(\mathbb{R}^n)$ 
is the classical inhomogeneous Besov or inhomogeneous Triebel--Lizorkin space.

Let  $f \in {\mathcal S}'$, then
we define  some sequences 
indexed by dyadic cubes $P$:

\bigskip

$c(\tilde{B}^{s'}_{pq})^{\sigma}_{x_{0}}(P)=\\
(\sum_{i  \geq (-\log_{2}l(P))\vee 0}||2^{is'}
|f*\phi_{i}(x)|(2^{-i}+|x_{0}-x|)^{-\sigma}||_{L^p(P)}^{q})^{1/q}$, 
\\
$0 < p,\ q \leq \infty$, 

c($\tilde{F}^{s'}_{pq})^{\sigma}_{x_{0}}(P)=\\
||\{\sum_{i  \geq (-\log_{2}l(P))\vee 0}(2^{is'}
|f*\phi_{i}(x)|(2^{-i}+|x_{0}-x|)^{-\sigma})^q\}^{1/q}||_{L^p(P)} $, 
\\
 $0< p < \infty$,\ $0 <  q \leq \infty$, 
 
$c(\tilde{F}^{s'}_{\infty q})^{\sigma}_{x_{0}}(P)=\\
l(P)^{-\frac{n}{q}}
||\{\sum_{i \geq (-\log_{2}l(P))\vee 0}(2^{is'}
|f*\phi_{i}(x)|(2^{-i}+|x_{0}-x|)^{-\sigma})^q\}^{1/q}||_{L^q(P)}$,
\\ $0 <  q \leq \infty$,

\noindent
with  the usual modification for $q=\infty$.

We  shall use the notation 
 $\tilde{E}^{s'}_{pq}$ with either $\tilde{B}^{s'}_{pq}$\ 
 or $\tilde{F}^{s'}_{pq}$. 
 We say the B-type case when $\tilde{E}^{s'}_{pq}=\tilde{B}^{s'}_{pq}$, and the F-type case when $\tilde{E}^{s'}_{pq}=\tilde{F}^{s'}_{pq}$.

\bigskip

\noindent
{\bf Definition 2.}\ \ 
Let $s, \ s', \ \sigma\ \in {\Bbb R},\  0 < p,q \leq \infty$ and 
  $x_{0} \in \mathbb{R}^n$.

The space $A^s(\tilde{E}^{s'}_{pq})^{\sigma}_{x_{0}}$ is defined to be the space of all 
$f \in \mathcal{S}'$ such that 
$$ 
||f||_{A^s(\tilde{E}^{s'}_{pq})^{\sigma}_{x_{0}}}\equiv
\sup_{\mathcal{D}\ni P}l(P)^{-s}c(\tilde{E}^{s'}_{pq})^{\sigma}_{x_{0}}(P) < \infty.
$$
The space $A^s(\tilde{E}^{s'}_{pq})^{\sigma}_{x_{0}}$ is the classical 2--microlocal Besov or  Triebel--Lizorkin space.

\noindent
We use the abbreviation 
$A^0(\tilde{{E}}^{s'}_{pq})^{\sigma}_{x_{0}}\equiv 
(\tilde{{E}}^{s'}_{pq})^{\sigma}_{x_{0}}$.  

\bigskip

\noindent
{\bf Examples.}

\begin{enumerate}[(i)]
\item \  The spaces 
$A^0(E^{s'}_{pq})^{0}_{x_{0}}= A^0(\tilde{E}^{s'}_{pq})^{0}_{x_{0}}= E^{s'}_{pq}(\mathbb{R}^n)$ are  the inhomogeneous Besov spaces or  
inhomogeneous Triebel--Lizorkin spaces [37], [38].

\item \  
The Besov type spaces  $B^{s, \tau}_{pq}(\mathbb{R}^n)$ and 
the Triebel--Lizorkin type spaces  $F^{s, \tau}_{pq}(\mathbb{R}^n)$ 
introduced by  D. Yang , W. Sickel and W. Yuan [36] ,
are contained in our definition as  

\noindent
$ E^{s', s}_{pq}(\mathbb{R}^n) =A^{ns}(E^{s'}_{pq})^{0}_{x_{0}}=A^{ns}(\tilde{E}^{s'}_{pq})^{0}_{x_{0}}$.

\item \  The Besov-Morrey spaces $\mathcal{N}^s_{uqp}$, and the Triebel--Lizorkin-Morrey spaces 
 $\mathcal{E}^s_{uqp}$  
studied by  Y. Sawano and H. Tanaka [34], or Y. Sawano, D. Yang and W. Yuan [35]
are realized in our definition as

 $\mathcal{N}^s_{uqp}\subset A^{n(\frac{1}{p}-\frac{1}{u})}
(B^s_{pq})^0_{x_{0}}$ if $0 < p \leq u \leq \infty$ and $0 < q \leq \infty$, 

 $\mathcal{E}^s_{uqp}=A^{n(\frac{1}{p}-\frac{1}{u})}
(F^s_{pq})^0_{x_{0}}$ if $0 < p \leq u \leq \infty$ and $0 < q \leq \infty$. 

The Morrey space ${\mathcal M}^u_{p}$ is realized as

${\mathcal M}^u_{p}= 
 A^{n(\frac{1}{p}-\frac{1}{u})}(F^0_{p2})^0_{x_{0}}$ if 
$1 < p < u < \infty$.

\item \  The $\dot{B}_{\sigma}$-Morrey spaces $\dot{B}_{\sigma}(L_{p,\lambda})$
studied by Y. Komori-Furuya et al. [28], 
are contained in our definition as  
 
$\dot{B}_{\sigma}(L_{p,\lambda})=
A^{\lambda+\frac{n}{p}}(F^0_{p2})^{\sigma}_{0}$, \ $1 < p < \infty$.

\item \  The 2-microlocal Besov spaces 
$B^{s, s'}_{pq}(U)$ studied  in H. Kempka [23, 24], 
are realized in our definition  as 

$B^{s, s'}_{pq}(U)=(\tilde{B}^{s+s'}_{pq})^{-s'}_{x_{0}}$  
 when $U=\{ x_{0} \}$.

\item \ The local Morrey spaces   
 $LM_{p, \lambda}$ introduced by V.I. Burenkov and H.V. Guliyes [6] and 
studied in Ts. Batbold and Y. Sawano [2], 
are realized in our definition as 

$LM_{p, \lambda}=(F^{0}_{p2})^{\lambda/p}_{0}$, \ $1 < p < \infty$.

\item \ The spaces $C^{s,s'}_{x_{0}}$ studied in  Y. Meyer [31], [32], 
are realized in our definition as 

$C^{s,s'}_{x_{0}}=
(\tilde{B}^{s+s'}_{\infty \infty})^{-s'}_{x_{0}}$=
$(B^{s+s'}_{\infty \infty})^{-s'}_{x_{0}}$.

\end{enumerate}

\bigskip

\section{Sequence spaces}

\noindent
For a sequence $c=(c(R))$ with $l(R) \leq 1$ we define  some sequences indexed by dyadic cubes $P$:

\bigskip

$
c(b^s_{pq})(P)
=(\sum_{i \geq (-\log_{2}l(P))\vee 0}||\sum_{l(R)=2^{-i}}2^{is}
|c(R)|\chi_{R}||_{L^p(P)}^{q})^{1/q}$, 
\\
$ 0 < p,\ q \leq \infty$,

$
c(f^s_{pq})(P)=||\bigl\{\sum_{i \geq (-\log_{2}l(P))\vee 0}\bigl(
\sum_{l(R)=2^{-i}}2^{is}
|c(R)|\chi_{R}\bigr)^q\bigr\}^{1/q}||_{L^p(P)}$, 
\\
$ 0 < p < \infty,\ 0 < q \leq \infty$,
and

$
c(f^s_{\infty q})(P)=l(P)^{-\frac{n}{q}}
||\bigl\{\sum_{i \geq (-\log_{2}l(P))\vee 0}\bigl(\sum_{l(R)=2^{-i}}2^{is}
|c(R)|\chi_{R}\bigr)^q\bigr\}^{1/q}||_{L^q(P)}$,

\noindent
$ 0 < q \leq \infty$, with the usual modification for $q=\infty$.

 The notation
 $e^s_{pq}$ is used to denote either $b^s_{pq}$ or 
$f^s_{pq}$. 
We say the B-type case when $e^{s'}_{pq}=b^{s'}_{pq}$, and the F-type case when $e^{s'}_{pq}=f^{s'}_{pq}$.
\bigskip

\noindent
{\bf Definition 3.}
\ \ Let $s, \ s' ,\ \sigma\ \in {\Bbb R},\ 0< p,q \leq \infty$ and 
 $x_{0} \in \mathbb{R}^n$.

We define the sequence space $a^s(e^{s'}_{pq})^{\sigma}_{x_{0}}$ 
to be the space of all sequences  $c=(c(R))_{l(R)\leq 1}$ such that 
$$
||c||_{a^s(e^{s'}_{pq})^{\sigma}_{x_{0}}}\equiv 
\sup_{\mathcal{D}\ni Q\ni x_{0} }l(Q)^{-\sigma}
\sup_{\mathcal{D} \ni P \subset 3Q\ }l(P)^{-s}c(e^{s'}_{pq})(P) < \infty.
$$

\bigskip

We use the abbreviation 
$a^0(e^{s'}_{pq})^{\sigma}_{x_{0}} \equiv 
(e^{s'}_{pq})^{\sigma}_{x_{0}}$, $a^s(e^{s'}_{pq})^{0}_{x_{0}}\equiv a^s(e^{s'}_{pq})$ and

\noindent
 $a^0(e^{s'}_{pq})^{0}_{x_{0}}\equiv e^{s'}_{pq}
\equiv e^{s'}_{pq}(\mathbb{R}^n)$. 
We note that the space 
$a^s(e^{s'}_{pq})$ is the sequence space of the inhomogeneous Besov type space or 
the inhomogeneous Triebel--Lizorkin type space  in the sense of Yang--Sickel--Yuan [36] and  the space $e^{s'}_{pq}
\equiv e^{s'}_{pq}(\mathbb{R}^n)$ 
is the sequence space of the classical inhomogeneous Besov or inhomogeneous Triebel--Lizorkin space.

\bigskip

\noindent
{\bf Remark 1.}\  It is easy that
when $\sigma < 0$, we have 
\ \ 
$A^s(E^{s'}_{pq})^{\sigma}_{x_{0}}= \{ 0 \}$ and\ \ $a^s(e^{s'}_{pq})^{\sigma}_{x_{0}}= \{ 0 \}$\ 
for\ \ $0 < p, q \leq \infty$ (See Proposition 1 below).

\bigskip


We define that for a sequence $(c(R))_{l(R)\leq 1}$,

\bigskip

$c(\tilde{b}^{s'}_{pq})^{\sigma}_{x_{0}}(P)=\\
(\sum_{i \geq (-\log_{2}l(P))\vee 0}||\sum_{l(R)=2^{-i}}2^{is'}|c(R)|
(2^{-i}+|x_{0}-x|)^{-\sigma}\chi_{R}||_{L^p(P)}^{q})^{1/q}$, 
\\
$0 < p,\ q \leq \infty$, 

$c(\tilde{f}^{s'}_{pq})^{\sigma}_{x_{0}}(P)=\\
||\{\sum_{i \geq (-\log_{2}l(P))\vee 0}(\sum_{l(R)=2^{-i}}2^{is'}
|c(R)|(2^{-i}+|x_{0}-x|)^{-\sigma}\chi_{R})^q\}^{1/q}||_{L^p(P)} $, 
\\
 $0< p < \infty$,\ $0 <  q \leq \infty$,

$c(\tilde{f}^{s'}_{\infty q})^{\sigma}_{x_{0}}(P)=
l(P)^{-\frac{n}{q}}\times\\
||\{\sum_{i \geq (-\log_{2}l(P))\vee 0}(\sum_{l(R)=2^{-i}}2^{is'}|c(R)|
(2^{-i}+|x_{0}-x|)^{-\sigma}\chi_{R})^q\}^{1/q}||_{L^q(P)}$,
\\
$0 <  q \leq \infty$,

 \noindent
 with the usual modification for $q=\infty$.

The notation
 $\tilde{e}^{s'}_{pq}$ is used  to denote either $\tilde{b}^{s'}_{pq}$\ 
 or $\tilde{f}^{s'}_{pq}$. 
We say the B-type case when $\tilde{e}^{s'}_{pq}=\tilde{b}^{s'}_{pq}$, and the F-type case when $\tilde{e}^{s'}_{pq}=\tilde{f}^{s'}_{pq}$.

\bigskip

\noindent
{\bf Definition 4.}
\ \ Let $s, \ s' ,\ \sigma\ \in {\Bbb R},\ 0< p,q \leq \infty$ and 
 $x_{0} \in \mathbb{R}^n$.

We define the sequence space  
$a^s(\tilde{e}^{s'}_{pq})^{\sigma}_{x_{0}}$ to be the space of all sequences 
$ c=(c(R))_{l(R) \leq 1}$ such that
$$ 
||c||_{a^s(\tilde{e}^{s'}_{pq})^{\sigma}_{x_{0}}}\equiv 
\sup_{\mathcal{D} \ni P}l(P)^{-s}c(\tilde{e}^{s'}_{pq})^{\sigma}_{x_{0}}(P) < 
\infty. 
$$

\noindent
We use the abbreviation 
$a^0(\tilde{{e}}^{s'}_{pq})^{\sigma}_{x_{0}}\equiv 
(\tilde{{e}}^{s'}_{pq})^{\sigma}_{x_{0}}$.  

\bigskip

\noindent
{\bf Definition 5.}\ \ 
Let $r_{1},\ r_{2}  \geq 0$  and $L > 0$. 
We say that a matrix operator $A=\{ a_{QP} \}_{Q P}$, 
indexed by dyadic cubes $Q$ and $P$, 
is  ($r_{1}, r_{2}, L$)-almost diagonal if the matrix $\{ a_{QP} \}$ 
satisfies  

\bigskip

$|a_{QP}| \leq C\bigl( \frac{l(Q)}{l(P)} \bigl)^{r_{1}}(1+l(P)^{-1}|x_{Q}-x_{P}|)^{-L}$ if $l(Q) \leq l(P)$,

 $
|a_{Q P}| \leq C\bigl( \frac{l(P)}{l(Q)} \bigl)^{r_{2}}(1+l(Q)^{-1}|x_{Q}-x_{P}|)^{-L}$ if $l(Q) > l(P)$.

\bigskip

The results about the boundedness of almost diagonal operators in 
 [9: Theorem 3.3], also
hold in our cases.
  
\bigskip

\noindent
{\bf Lemma\ 1}.\ \ {\it 
Suppose that $ s,\ s',\ \sigma\ \in {\Bbb R}, \ \ x_{0} \in \mathbb{R}^n$ 
and 
\ $0 < p,\ q\ \leq \infty$. Then,

{\rm (i)}\ \  an $(r_{1}, r_{2}, L)$--almost diagonal 
matrix operator A is bounded on  $a^s(e^{s'}_{pq})^{\sigma}_{x_{0}}$ 
for $r_{1} > \max(s',\ \sigma+s + s'- \frac{n}{p}),\ r_{2} > J-s'$ and $L>J$ 
where 
$J= n/ \min(1,\ p, \ q)$ 
in the case $e^{s'}_{pq}=f^{s'}_{pq}$ , and 
$J= n/ \min(1,\ p)$ in the case $e^{s'}_{pq}=b^{s'}_{pq}$, respectively,

{\rm (ii)}\ \ 
 an $(r_{1}, r_{2}, L)$-almost diagonal 
matrix operator A  is bounded on 
$a^s(\tilde{e}^{s'}_{pq})^{\sigma}_{x_{0}}$ 
for $r_{1} > \max(s'+(\sigma\vee 0),\ \ (\sigma\vee 0)+s + s'- \frac{n}{p}),\ r_{2} > J-s'+(\sigma \wedge 0)$ and $L>J$ 
where $J= n/ \min(1,\ p, \ q)$ 
in the case  $\tilde{e}^{s'}_{pq}=\tilde{f}^{s'}_{pq}$, and 
$J= n/ \min(1,\ p)$ in the case 
$\tilde{e}^{s'}_{pq}=\tilde{b}^{s'}_{pq}$, respectively.

}

\bigskip

\noindent
{\it Proof} :\ \ (i)\ \ We may assume $\sigma \geq 0$ by Remark 1.
We assume that $A = (a_{R R'})$ is $(r_{1}, r_{2}, L)$-- almost diagonal. 
Let $c= (c(R)) \in a^s(e^{s'}_{pq})^{\sigma}_{x_{0}}$. 
For  dyadic cubes $P$ and $R$ with $R \subset P$, 
we write $Ac= A_{0}c + A_{1}c + A_{2}c$ with 

\begin{eqnarray*}
&&
(A_{0}c)(R)=\sum_{l(R) \leq l(R') \leq l(P)}a_{R R'}c(R'), 
\\
&&
(A_{1}c)(R)=\sum_{l(R')  < l(R) \leq l(P)}a_{R R'}c(R'), 
\\
&&
(A_{2}c)(R)=\sum_{l(R) \leq l(P) < l(R') \leq 1 }a_{R R'}c(R').
\end{eqnarray*}
We claim that
$$
||A_{i}c||_{a^{s}(e^{s'}_{pq})^{\sigma}_{x_{0}}}
\leq 
C||c||_{a^{s}(e^{s'}_{pq})^{\sigma}_{x_{0}}}, \ \ i=0,1,2.
$$
We will consider the case of F-type  for $0 < p < \infty$, $0 < q \leq \infty$.
Since $A$ is almost diagonal, we see that 
for dyadic cubes $P$ with $l(P)=2^{-j}$, 
\begin{eqnarray*}
\lefteqn{(A_{0}c)(f^{s'}_{pq})(P) = 
||\bigl\{\sum_{i \geq j\vee 0}\sum_{l(R)=2^{-i}}\bigl(2^{is'}|(A_{0}c)(R)|
\bigr)^{q}\chi_{R}\bigr\}^{1/q}||_{L^{p}(P)}
} 
\\
&\leq& 
C||\bigl\{\sum_{i \geq j\vee 0}\sum_{l(R)=2^{-i}}2^{is'q}\bigl(
\sum_{i \geq k \geq j\vee 0}\sum_{l(R')=2^{-k}}|a_{R R'}||c(R')|\bigr)^{q}
\chi_{R}\bigr\}^{1/q}||_{L^{p}(P)} 
\\
&\leq&
C||\bigl\{\sum_{i \geq j\vee 0}\sum_{l(R)=2^{-i}}2^{is'q}\times 
\\
\lefteqn{\bigl(\sum_{i \geq k \geq j\vee 0}
\sum_{l(R')=2^{-k}}2^{-(i-k)r_{1}}(1+2^k|x_{R}-x_{R'}|)^{-L}
|c(R')|\bigr)^{q}\chi_{R}\bigr\}^{1/q}||_{L^{p}(P)}.
}
\end{eqnarray*}
Using the maximal function $M_{t}f(x)$,\ $0< t\leq 1$,  defined by 
$$
M_{t}f(x)=\sup_{x \in Q}\bigr(\frac{1}{l(Q)^n}\int_{Q}
|f(y)|^t\ dy\bigl)^{1/t}
$$
(cf. [28: Lemma 7.1] or [9: Remark A.3]), 
we have for $L > n/t$,
\begin{eqnarray*}
\lefteqn{
(A_{0}c)(f^{s'}_{pq})(P) \leq 
C||\bigl\{\sum_{i \geq j\vee 0}\sum_{l(R)=2^{-i}}2^{is'q}2^{-ir_{1}q}\times
}
\\
&&
\Bigl(\sum_{i \geq k \geq j\vee 0}2^{k r_{1}}2^{(k-i)_{+}n/t}M_{t}\bigl(\sum_{l(R')=2^{-k}}|c(R')|
\chi_{R'}\bigr)
\Bigr)^{q}\chi_{R}\bigr\}^{1/q}||_{L^{p}(P)} 
\\
\lefteqn{
\leq 
C||\big\{\sum_{i \geq j\vee 0}2^{-i(r_{1}-s')q}\Bigl(\sum_{i \geq k \geq j\vee 0}
2^{kr_{1}}
M_{t}\bigl(\sum_{l(R')=2^{-k}}|c(R')|\chi_{R'}\bigr)\Bigr)
^{q}\bigr\}^{1/q}||
_{L^{p}(P)} 
}
\\
&\leq& 
C||\bigl\{\sum_{i \geq j\vee 0}2^{is'q}M_{t}\bigl(\sum_{l(R')=2^{-i}}|c(R')|
\chi_{R'}\bigr)^{q}
\bigr\}^{1/q}||_{L^{p}(P)} 
\\
&\leq &
C||\big\{\sum_{i \geq j\vee 0}2^{is'q}\bigl(\sum_{l(R')=2^{-i}}|c(R')|\chi_{R'}
\bigr)^{q}
\bigr\}^{1/q}||_{L^{p}(P)} 
=Cc(f^{s'}_{pq})(P),
\end{eqnarray*}
where these inequalities follow from Hardy's inequality if $r_{1} > s'$ and the  Fefferman-Stein inequality if $0< t< \min(p,q)$.

For the B-type case  we have the same estimate for $r_{1} > s'$ and $0 < t< \min(1,p)$.

Therefore, we get the estimate 

$$
A_{0}c(e^{s'}_{pq})(P)
\leq 
Cc(e^{s'}_{pq})(P)
$$
if $r_{1} > s'$,  $0 < p < \infty$, $0 < q \leq \infty$, $L > J$.

In the sane way  we will get the estimate for $(A_{1}c)(f^{s'}_{pq})(P)$.
 We have that for dyadic cubes $P$ with $l(P)=2^{-j}$, 
\begin{eqnarray*}
\lefteqn{(A_{1}c)(f^{s'}_{pq})(P) = 
||\bigl\{\sum_{i \geq j\vee 0}\sum_{l(R)=2^{-i}}\bigl(2^{is'}|(A_{1}c)(R)|
\bigr)^{q}\chi_{R}\bigr\}^{1/q}||_{L^{p}(P)}
} 
\\
&\leq &
C||\bigl\{\sum_{i \geq j\vee 0}\sum_{l(R)=2^{-i}}2^{is'q}\bigl(
\sum_{i \leq k }\sum_{l(R')=2^{-k}}|a_{R R'}||c(R')|\bigr)^{q}
\chi_{R}\bigr\}^{1/q}||_{L^{p}(P)} 
\\
&\leq&
C||\bigl\{\sum_{i \geq j\vee 0}\sum_{l(R)=2^{-i}}2^{is'q}\times 
\\
\lefteqn{\bigl(\sum_{i \leq k }
\sum_{l(R')=2^{-k}}2^{-(k-i)r_{2}}(1+2^i|x_{R}-x_{R'}|)^{-L}
|c(R')|\bigr)^{q}\chi_{R}\bigr\}^{1/q}||_{L^{p}(P)}.
}
\end{eqnarray*}
Using the maximal function $M_{t}f(x)$ as  above, 
we have 
\begin{eqnarray*}
\lefteqn{
(A_{1}c)(f^{s'}_{pq})(P) \leq 
C||\bigl\{\sum_{i \geq j\vee 0}\sum_{l(R)=2^{-i}}2^{is'q}2^{ir_{2}q}\times
}
\\
&&
\Bigl(\sum_{i \leq k }2^{-k r_{2}}2^{(k-i)_{+}n/t}M_{t}\bigl(\sum_{l(R')=2^{-k}}|c(R')|
\chi_{R'}\bigr)
\Bigr)^{q}\chi_{R}\bigr\}^{1/q}||_{L^{p}(P)} 
\\
&\leq& 
C||\big\{\sum_{i \geq j\vee 0}2^{i(r_{2}+s'-n/t)q}\times
\\
&&
\Bigl(\sum_{i \leq k }
2^{-k(r_{2}-n/t)}M_{t}\bigl(\sum_{l(R')=2^{-k}}|c(R')|\chi_{R'}\bigr)\Bigr)
^{q}\bigr\}^{1/q}||
_{L^{p}(P)} 
\\
&\leq& 
C||\bigl\{\sum_{i \geq j\vee 0}2^{is'q}M_{t}\bigl(\sum_{l(R')=2^{-i}}|c(R')|
\chi_{R'}\bigr)^{q}
\bigr\}^{1/q}||_{L^{p}(P)} 
\\
&\leq& 
C||\big\{\sum_{i \geq j\vee 0}2^{is'q}\bigl(\sum_{l(R')=2^{-i}}|c(R')|\chi_{R'}
\bigr)^{q}
\bigr\}^{1/q}||_{L^{p}(P)} 
=Cc(f^{s'}_{pq})(P),
\end{eqnarray*}
where these inequalities follow from Hardy's inequality if $r_{2} + s'-n/t > 0$
and the  Fefferman-Stein inequality if $0 < t < \min(p,q)$.

In the same way we get the same estimate for the B-type case   that 
$$
(A_{1}c)(b^{s'}_{pq})(P) \leq Cc(b^{s'}_{pq})(P)
$$
if $r_{2} +s'- n/t>0$, $0 < t < \min(1, p)$. 
Therefore, we get the estimate 

$$
A_{1}c(e^{s'}_{pq})(P)
\leq 
Cc(e^{s'}_{pq})(P)
$$
if $r_{2} >J- s'$, 
 $0 < p < \infty$, $0 < q \leq \infty$, $L > J$.

When $p=\infty$, we get the same estimate.
Thus, we get 

$$
||A_{i}c||_{a^{s}(e^{s'}_{pq})^{\sigma}_{x_{0}}}
\leq 
C||c||_{a^{s}(e^{s'}_{pq})^{\sigma}_{x_{0}}}, \ \ i=0,1
$$

\noindent
if $r_{1} > s'$, $r_{2} >J-s'$, $L>J$, $0 < p \leq \infty$ and $0 < q \leq \infty$.
 
Next, we will give the estimates for the $A_{2}$ case.

\noindent
We note that  if $L > n$ , 
$$\sum_{l(P)=2^{-j}}(1+2^{j}|x_{R}-x_{P}|)^{-L} < \infty$$ 
(cf. [4, Lemma 3.4]), and 
if $c \in a^{s}(e^{s'}_{pq})^{\sigma}_{x_{0}}$, then 

$$
|c(R)| \leq C(|x_{0}-x_{R}|+l(R))^{\sigma}l(R)^{s+s'-n/p}
||c||_{a^{s}(e^{s'}_{pq})^{\sigma}_{x_{0}}}
$$
for a dyadic cube $R \subset 3 Q$  and $x_{0} \in Q$. 
Hence, we obtain, for dyadic cubes $P$ with $l(P)=2^{-j}$, $0 < p <\infty$ and $0 < q \leq \infty$, 
\begin{eqnarray*}
\lefteqn{
(A_{2}c)(f^{s'}_{pq})(P) =
||\bigl\{\sum_{i \geq j}\sum_{l(R)=2^{-i}}\bigl(2^{is'}|(A_{2}c)(R)|
\bigr)^{q}
\chi_{R}\bigr\}^{1/q}||_{L^{p}(P)} 
}
\\
\lefteqn{
\leq 
C||\bigl\{\sum_{i \geq j}\sum_{l(R)=2^{-i}}2^{is'q}\times
}
\\
\lefteqn{
\bigl(\sum_{j \geq k \geq 0}\sum_{l(R')=2^{-k}}2^{-(i-k)r_{1}}
(1+2^{k}|x_{R}-x_{R'}|
\bigr)^{-L}|c(R')|)^{q}\chi_{R}\bigr\}^{1/q}||_{L^{p}(P)} 
}
\\
\lefteqn{
\leq
C||\bigl\{\sum_{i \geq j}2^{-i(r_{1}-s')q}\times 
}
\\
&&
\bigl(
\sum_{j \geq k \geq 0}
2^{kr_{1}}
2^{-k(\sigma + s + s'-n/p)}(1+2^k|x_{0}-x_{P}|)^{\sigma}||c||_{a^{s}(f^{s'}_{pq})^{\sigma}_{x_{0}}}\bigr)^{q}\bigr\}^{1/q}||_{L^{p}(P)}
\\
\lefteqn{\leq
C2^{-j(r_{1}-s')}2^{-jn/p}
\sum_{j\geq k \geq 0}2^{k(r_{1}-\sigma- s -s'+ n/p)}
(1+2^{j}|x_{0}-x_{P}|)^{\sigma}
||c||_{a^{s}(f^{s'}_{pq})^{\sigma}_{x_{0}}}
}
\\
\lefteqn{\leq
C2^{-j(r_{1}-s'+n/p)}2^{j(r_{1}-\sigma- s -s'+ n/p)}(1+2^{j}|x_{0}-x_{P}|)^{\sigma}||c||_{a^{s}(f^{s'}_{pq})^{\sigma}_{x_{0}}}
}
\\
\lefteqn{\leq 
C2^{-js}(2^{-j}+|x_{0}-x_{P}|)^{\sigma}||c||_{a^{s}(f^{s'}_{pq})^{\sigma}_{x_{0}}} 
}
\end{eqnarray*}
where these inequalities follow 
if $r_{1} >\sigma +s+s'-\frac{n}{p}$,  $r_{1} > s'$, $ L > n$
and $\sigma \geq 0$.

In the same way for the B-type case  we have the same estimate.

Hence, we have,  
 
$$
||A_{2}c||_{a^{s}(e^{s'}_{pq})^{\sigma}_{x_{0}}} 
\leq
C||c||_{a^{s}(e^{s'}_{pq})^{\sigma}_{x_{0}}}
$$

\noindent
if $r_{1} > \sigma + s+s'-n/p$, $r_{1} > s'$,  $0 < p < \infty$ and $0 < q \leq \infty$.  

We get the same estimate for the case $p = \infty$. 
Therefore, we obtain the desired conclusion.

(ii)\ \ We put $w_{i} =(2^{-i}+|x_{0}-x|)^{-\sigma}$. 
We see that 
 $w_{i}\leq 2^{(i-k)_{+}\sigma}w_{k}$ if  $0 \leq \sigma$, and 
 $w_{i}\leq 2^{(k-i)_{+}\sigma}w_{k}$ if  $0 > \sigma$. Then, using these inequalities  
we can prove  the desired result by using the same way in the above proof of (i).
\qed

\bigskip

\noindent
{\bf Lemma 2.}\ \ {\it 
 Let $r_{1}, r_{2} \in \mathbb{N}_{0}, \ L >n $ and $L_{1} > n+r_{1}, L_{2} > n + r_{2}$. 
 Assume that  for dyadic cubes $P$ and $R$, 
$\phi_{P}$ and $\varphi_{R}$ 
 are functions on $\mathbb{R}^{n}$ satisfying following properties: 
\begin{eqnarray*}
&&
(2.1) \ \ \int_{\mathbb{R}^n}\phi_{P}(x)x^{\gamma} dx =0\ \  
 for\  \ |\gamma| < r_{1},
\\
&&
(2.2)\ \ |\phi_{P}(x)| \leq C(1+l(P)^{-1}|x-x_{P}|)^{-\max(L, L_{1})},
\\
&&
(2.3)\ \ |\partial^{\gamma}\phi_{P}(x)| \leq Cl(P)^{-|\gamma|}(1+l(P)^{-1}|x-x_{P}|)^{-L}\ \ 
\\
&& for\ \  0<|\gamma| \leq r_{2},
\\
&&
(2.4)\ \ \int_{\mathbb{R}^n}\varphi_{R}(x) x^{\gamma} dx =0\ \ 
for \ \ |\gamma|< r_{2}, 
\\
&&
(2.5)\ \ |\varphi_{R}(x)| \leq C (1+l(R)^{-1}|x-x_{R}|)
^{-\max(L, L_{2})},
\\
&&
(2.6)\ \ |\partial^{\gamma}\varphi_{R}(x)| \leq C l(R)^{-|\gamma|}
(1+l(R)^{-1}|x-x_{R}|)^{-L} \ \ 
\\
&& for\ \  0<|\gamma| \leq r_{1},
\end{eqnarray*}
where  {\rm (2.1)} and  {\rm (2.6)} are void when $r_{1}=0$, and 
{\rm (2.3)} and {\rm (2.4)} are void when  $r_{2}=0$. 
Then, we have that

\bigskip

$ l(P)^{-n}|
\langle \phi_{P}\ ,\ \varphi_{R} \rangle| 
\leq C \bigl( \frac{l(P)}{l(R)} \bigl)^{r_{1}}(1+l(R)^{-1}|x_{P}-x_{R}|)^{-L}$

  if \ \ $l(P) \leq l(R)$, 

 $ l(R)^{-n}|\langle \phi_{P}\ ,\ \varphi_{R} \rangle| 
\leq C \bigl( \frac{l(R)}{l(P)} \bigl)^{r_{2}}(1+l(P)^{-1}|x_{P}-x_{R}|)^{-L}$

if\ \  $l(R) < l(P)$.
}
\bigskip

\noindent
{\it Proof}.\ \  We refer to [10: Corollary  B.3] , [5: Lemma 6.3] or  [29: Lemma 3.1].
\qed

\bigskip

\noindent
{\bf Lemma\ 3}.\ \ {\it 
Suppose that $ s,\ s', \ \sigma \in {\Bbb R}, \ x_{0} \in \mathbb{R}^n$ 
and 
\ $0 < p,\ q\ \leq \infty$. Let $r_{1}, r_{2} \in \mathbb{N}_{0}$ and $L >n$. 
Assume that  functions $\phi_{P}$ and $\varphi_{P}$ satisfy 
 {\rm (2.1)}, {\rm (2.2)}, {\rm (2.3)}, {\rm (2.4)}, {\rm (2.5)}, {\rm (2.6)}  in Lemma {\rm 2}. Let $J$  as in Lemma {\rm 1}.
Then we have

${\rm (i)}$\ \ for  a dyadic cube $R$ and a sequence 
$c \in a^s(e^{s'}_{pq})^{\sigma}_{x_{0}}$, 

$\sum_{\mathcal{D} \ni P,\  l(P) \leq 1}c(P)\langle \phi_{P}\ ,\ \varphi_{R}\rangle$ is convergent  if $r_{1} > J-n-s'$ and $L>J$,

${\rm (ii)}$\ \ for  a dyadic cube $R$ and a sequence $c \in a^s(\tilde{e}^{s'}_{pq})^{\sigma}_{x_{0}}$,

$\sum_{\mathcal{D} \ni P,\ l(P) \leq 1}c(P) \langle \phi_{P}\ ,\ \varphi_{R}\rangle$ is convergent  if $r_{1} > J-n-s'- (\sigma \wedge 0)$ and 
$L>J+\sigma$.
}

\bigskip

\noindent
{\it Proof} :\ \ 
(i)\ \ We may assume that \ $\sigma \geq 0$ by Remark 1. 

\noindent
We write  
$\sum_{\mathcal{D} \ni P}c(P)\langle \phi_{P}\ ,\ \varphi_{R} \rangle
=I=  I_{0} + I_{1}$ with 

\begin{eqnarray*}
&&
I_{0}=\sum_{l(R) \leq l(P) \leq 1}c(P)\langle \phi_{P}\ ,\ \varphi_{R} \rangle, 
\\
&&
I_{1}=\sum_{l(P)  < l(R) }c(P)\langle \phi_{P}\ ,\ \varphi_{R} \rangle 
\end{eqnarray*}
for $c \in a^s(e^{s'}_{pq})^{\sigma}_{x_{0}}$. We claim that
$
I_{i}< \infty, \ \ i=0,1.
$

For a dyadic cube $R$ with $l(R)=2^{-i}$ we have, by Lemma 2 that
\begin{eqnarray*}
\lefteqn{
|I_{0}| \leq C
\sum_{i \geq j\geq 0}\sum_{l(P)=2^{-j}}|c(P)||\langle \phi_{P}\ ,\ \varphi_{R} \rangle|
}
\\
&\leq& 
C\sum_{i \geq j\geq 0}\sum_{l(P)=2^{-j}}|c(P)
|2^{-in}2^{(j-i)r_{2}}(1+2^{j}|x_{R}-x_{P}|)^{-L}
\\
&\leq&
C\sum_{i \geq j\geq 0}2^{-i(r_{2}+n)}
2^{jr_{2}}M_{t}(\sum_{l(P)=2^{-j}}|c(P)|\chi_{P})(x),
\end{eqnarray*}
for $L > n/t$, $0< t< 1$ and $x \in R$.
Taking $L^1(R)$ norm and using the  Fefferman-Stein inequality, 
we have, 
\begin{eqnarray*}
\lefteqn{
|I_{0}|2^{-in} = ||I_{0}||_{L^1(R)} 
}
\\
&\leq& C 2^{-in}
||\sum_{i \geq j\geq 0}M_{t}(\sum_{l(P)=2^{-j}}|c(P)|\chi_{P})||
_{L^{1}(R)}
\\
&\leq& 
C2^{-in}||
\sum_{i \geq j \geq 0}\sum_{l(P)=2^{-j}}|c(P)|
\chi_{P}||_{L^{1}(R)}
\\
&\leq&
C\sum_{1 \geq l(P),\ R \subset P }|c(P)|2^{-2in}< \infty.
\end{eqnarray*}

\bigskip

In the same way we obtain the estimate of $I_{1}$:
\begin{eqnarray*}
\lefteqn{
|I_{1}| \leq C
\sum_{j \geq i\vee 0}\sum_{l(P)=2^{-j}}|c(P)||\langle \phi_{P}\ ,\ \varphi_{R} \rangle|
}
\\
&\leq& 
C\sum_{j \geq i\vee 0}\sum_{l(P)=2^{-j}}|c(P)
|2^{-jn}2^{(i-j)r_{1}}(1+2^{i}|x_{R}-x_{P}|)^{-L}
\\
&\leq&
C\sum_{j \geq i\vee 0}2^{-j(r_{1}+n)}2^{ir_{1}}
\sum_{l(P)=2^{-j}}|c(P)|(1+2^{i}|x_{R}-x_{P}|)^{-L}
\\
&\leq&
C\sum_{j \geq i\vee 0}
2^{-j(r_{1}+n-n/t+s')}2^{ir_{1}}2^{-in/t}
M_{t}(\sum_{l(P)=2^{-j}}2^{js'}|c(P)|\chi_{P})(x)
\end{eqnarray*}
if $0 < t \leq 1, \ L > n/t$ and $x \in R$ with $l(R)= 2^{-i}$.

By using  the monotonicity of $l^q$-norm  and 
 H${\rm \ddot o}$lder's inequality,
we get the following result, 
$$
|I_{1}|\leq C2^{-i(n+s')}\{\sum_{j \geq i\vee 0}(M_{t}(\sum_{l(P)=2^{-j}}2^{js'}
|c(P)|\chi_{P})(x))^q \}^{1/q}
$$
if $r_{1}+n-n/t+s' > 0,\ 0 < q \leq \infty$ and $x \in R$.

Taking $L^p(R)$ norm and using the  Fefferman-Stein inequality, 
we have, for a dyadic cube $R$ with $l(R)=2^{-i}$ and $c \in 
a^{s}(f^{s'}_{pq})^{\sigma}_{x_{0}}$,
\begin{eqnarray*}
\lefteqn{
|I_{1}|2^{-in/p} = ||I_{1}||_{L^p(R)} 
\leq
C2^{-i(n+s')}c(f^{s'}_{pq})(R)
}
\\
&\leq&
 C2^{-i(n+s'+\sigma +s)}||c||_{a^{s}(f^{s'}_{pq})^{\sigma}_{x_{0}}} < \infty
\end{eqnarray*}
if $0 < t < \min(p,q)$, $0 < p < \infty,\ 0 < q \leq \infty$.
In the same way we  get the same estimate for the case $p=\infty$. 
Furthermore, we obtain the same estimate for the B-type case if  $0 < t < p$, 
$0 < p \leq \infty,\ 0 < q \leq \infty$.
Therefore, we obtain that $I_{1}$ is convergent
if  $r_{1} > J-n-s'$ and $L>J$.

(ii)\ \ Let $I_{0}$ and $I_{1}$ be as in the proof of (i). Then 
by arguing as in the proof of (i), we have 
$I_{0} < \infty$ for $L > n$. 
 We put $w_{j}(P)=(2^{-j}+|x_{P}-x_{0}|)^{-\sigma}$ for a dyadic cube  $P$ with $l(P)=2^{-j}$.

Note  that 
$$
|c(P)| \leq Cl(P)^{s+s'-n/p}w_{j}(P)^{-1}
||c||_{a^{s}(\tilde{e}^{s'}_{pq})^{\sigma}_{x_{0}}}$$ 
for $c \in a^{s}(\tilde{e}^{s'}_{pq})^{\sigma}_{x_{0}}$. 
We have, by Lemma 2 for  a dyadic cube $R$ with $l(R)=2^{-i}$ and $\sigma \geq 0$, 

\begin{eqnarray*}
\lefteqn{
|I_{1}| \leq C
\sum_{j \geq i\vee 0}\sum_{l(P)=2^{-j}}|c(P)||\langle \phi_{P}\ ,\ \varphi_{R} \rangle|
}
\\
&\leq& 
C\sum_{j \geq i\vee 0}\sum_{l(P)=2^{-j}}|c(P)
|2^{-jn}2^{(i-j)r_{1}}(1+2^{i}|x_{R}-x_{P}|)^{-L}
\\
&\leq& 
C\sum_{j \geq i\vee 0}\sum_{l(P)=2^{-j}}|c(P)
|2^{-jn}2^{(i-j)r_{1}}w_{j}(P)w_{j}(P)^{-1}\times
\\
&&
(1+2^{i}|x_{R}-x_{P}|)^{-L}
\\
&\leq&
C\sum_{j \geq i\vee 0}2^{-j(r_{1}+n)}2^{ir_{1}}2^{-i\sigma}
\sum_{l(P)=2^{-j}}|c(P)|w_{j}(P)\times
\\
&&(1+2^{i}|x_{R}-x_{P}|)^{-(L-\sigma)}
\\
&\leq&
C\sum_{j \geq i\vee 0}
2^{-j(r_{1}+n-n/t+s')}2^{i(r_{1}+\sigma-n/t)}\times
\\
&&
M_{t}(\sum_{l(P)=2^{-j}}2^{js'}w_{j}(P)|c(P)|\chi_{P})(x).
\end{eqnarray*}
By using the same way as in the proof of (i), we get
\begin{eqnarray*}
\lefteqn{
|I_{1}|2^{-ip/n}\leq C2^{-i(n+s'+\sigma)}c(\tilde{e}^{s'}_{pq})^{\sigma}_{x_{0}}(R)
}\\
&&
\leq C2^{-i(n+s'+\sigma +s)}||c||_{a^{s}(\tilde{e}^{s'}_{pq})^{\sigma}_{x_{0}}} < \infty
\end{eqnarray*}
if $r_{1} > J-n-s'$ and $L > \sigma+J$. We also obtain the same estimate for the case $\sigma < 0$. \qed

\bigskip

 For a sequence $c(P)$ with $l(P)=2^{-j}$, we define the sequence $c^*(P)$ by 
$$
c^*(P) = \sum_{l(R)=2^{-j}}|c(R)|(1+2^{j}|x_{P}-x_{R}|)^{-L}
$$ 

\noindent 
for  $L>J$ where $J$ is as in Lemma 1.

We define 
for $f \in \mathcal{S}'$, $\gamma \in \mathbb{N}_{0}$ and a dyadic cube $P$ with $l(P)=2^{-j}$, the sequence $\inf_{\gamma}(f)(P)$ and 
$t_{\gamma}(P)$ 
by 

$
\inf_{\gamma}(f)(P)= \max\{\inf _{R \ni y}|\phi_{j}*f (y)|: R \subset P, l(R)=2^{-(\gamma+j)} \},\ \ 
$

$$
t_{\gamma}(P)=\inf_{P \ni y}|\phi_{j-\gamma}*f(y)|.
$$

\bigskip

\noindent
{\bf Lemma 4.}\ \  {\it 
  For $ s', \ \sigma \in {\Bbb R}, \ x_{0} \in \mathbb{R}^n,\ 0 < p,q \leq \infty,\ f \in \mathcal{S}'$ and a dyadic cube $P$ with $l(P)=2^{-j}$, we have 

${\rm (i)}$\ \ $$
c(e^{s'}_{pq})(P) \sim c^*(e^{s'}_{pq})(P),\ \ 
c(\tilde{e}^{s'}_{pq})^{\sigma}_{x_{0}}(P) 
\sim c^*(\tilde{e}^{s'}_{pq})^{\sigma}_{x_{0}}(P),
$$

${\rm (ii)}$\ \ 

\bigskip
$\ \ \ \ \ \ \ \ \ \ \ \ \ \ \ \  
\inf_{\gamma}(f)(P)\chi_{P}\leq C2^{\gamma L}\sum_{R\subset P, l(R)
=2^{-(\gamma+j)}}t_{\gamma}^*(R)\chi_{R}.
$

\bigskip
\noindent
for $\gamma$ sufficient large.
}
\bigskip

\noindent
{\it Proof}. (i)
\ \ It suffices to prove 
$$
c^*(e^{s'}_{pq})(P) \leq C c(e^{s'}_{pq})(P)
$$
since $|c(P)| \leq c^*(P)$.

Using the  Fefferman-Stein inequality, we have 
\begin{eqnarray*}
\lefteqn{
c^*(f^{s'}_{pq})(P) = 
||\{ \sum_{i \geq j\vee 0}
(2^{is'}\sum_{l(R)=2^{-i}}|c^*(R)|\chi_{R})^{q} \}^{1/q}||_{L^p(P)}
}
\\
&\leq& 
C||\{ \sum_{i \geq j\vee 0}
(2^{is'}\sum_{l(R)=2^{-i}}\sum_{l(R')=2^{-i}}|c(R')|\times
\\
&&
(1+2^{i}|x_{R}-x_{R'}|)^{-L}\chi_{R})^{q} \}^{1/q}||_{L^p(P)}
\\
&\leq& 
C||\{ \sum_{i \geq j\vee 0}
(2^{is'}\sum_{l(R)=2^{-i}}M_{t}(\sum_{l(R')=2^{-i}}|c(R')|
\chi_{R'})\chi_{R})^{q} \}^{1/q}||_{L^p(P)}
\\
&\leq&
C||\{ \sum_{i \geq j\vee 0}
(\sum_{l(R')=2^{-i}}2^{is'}|c(R')|\chi_{R'})^{q} \}^{1/q}||_{L^p(P)}
=Cc(f^{s'}_{pq})(P)
\end{eqnarray*}
if $0< t < \min(p,q)$, $L > n/t$ and $0< p < \infty, 0< q \leq \infty$. 
Moreover, for the $p=\infty$ case,  we have the same result. For the B-type case , we obtain the same result by the same argument as above. 
We also obtain the same result for the other case.

(ii)\ \ Let $R_{0}$ and $R$  in $P$ be cubes with $l(R_{0})=l(R)=2^{-(\gamma+j)}$. It suffices to show 
$$
t_{\gamma}(R_{0}) \leq C2^{\gamma L}t_{\gamma}^*(R).
$$
Since 
$$
1 \leq 2^L2^{\gamma L}(1+ 2^{\gamma +j}|x_{R}-x_{R_{0}}|)^{-L}, 
$$
we have 
\begin{eqnarray*}
\lefteqn{
t_{\gamma}(R_{0}) \leq Ct_{\gamma }(R_{0})2^{\gamma L}
(1+2^{\gamma +j}|x_{R}-x_{R_{0}}|)^{-L}
}
\\
&\leq& 
C2^{\gamma L}\sum_{l(R')=2^{-(\gamma+j)}}t_{\gamma}(R')
(1+2^{\gamma +j}|x_{R}-x_{R'}|)^{-L}=C2^{\gamma L}t_{\gamma}^*(R).
\end{eqnarray*}\qed
\bigskip

\section{Characterizations}

 \bigskip

\noindent
{\bf Remark 2.}\  (See [11: (3.20)] ). \ \ 
Let $\phi_{0}$ be  a Schwartz function satisfying (1.1) and (1.2) and let 
$\phi$ be  a Schwartz function satisfying (1.3) and (1.4).
 Then there exist a Schwartz function $\varphi_{0}$ 
satisfying the same conditions (1.1) and (1.2) and a Schwartz function $\varphi$  satisfying the same conditions (1.3) and (1.4) such that 

\bigskip

$
\sum_{j \in {\Bbb N_{0}}}\hat{\varphi}_{j}(\xi)
\hat{\phi}_{j}(\xi)=1$ for any $\xi$
where $\varphi_{j}(x)=2^{jn}\varphi(2^jx), \ j \in \mathbb{N}$.

\bigskip

\noindent
Hence we have
the $\varphi$-transform [8; Lemma 2.1] for $f \in \mathcal{S}'$ such that 
$$
f=\sum_{l(Q)\leq 1}l(Q)^{-n}\langle f ,\ \varphi_{Q} \rangle \phi_{Q},
$$
where 
$\phi_{Q}(x)=
\phi(l(Q)^{-1}(x-x_{Q}))$ and $\varphi_{Q}(x)=
\varphi(l(Q)^{-1}(x-x_{Q}))$ 
for a dyadic cube $Q$ with $l(Q) < 1$, and 
$\phi_{Q}(x)=
\phi_{0}(l(Q)^{-1}(x-x_{Q}))$ and $\varphi_{Q}(x)=
\varphi_{0}(l(Q)^{-1}(x-x_{Q}))$ 
 for a dyadic cube $Q$ with $l(Q) = 1$.

\bigskip

\noindent
{\bf Theorem\ 1.}\  \ {\it 
For $s, \ s',\ \sigma \in {\Bbb R},\  0 < p,\ q \leq \infty,\ 
x_{0} \in \mathbb{R}^n$ and $\phi_{0},\ \phi \in \mathcal{S}$ 
as in Remark {\rm 2}, 
we have 

{\rm (i)}\ \ 

$$
A^s(E^{s'}_{pq})^{\sigma}_{x_{0}} = 
\{ f= \sum_{l(Q) \leq 1}c(Q)\phi_{Q} : 
\ \ (c(Q)) \in a^s(e^{s'}_{pq})^{\sigma}_{x_{0}} \},
$$  and 

{\rm (ii)}\ \ 

$$
A^s(\tilde{E}^{s'}_{pq})^{\sigma}_{x_{0}} = 
\{ f= \sum_{l(Q) \leq 1}c(Q)\phi_{Q} : 
\ \ (c(Q)) \in a^s(\tilde{e}^{s'}_{pq})^{\sigma}_{x_{0}} \}.
$$
}

\bigskip

\noindent
{\bf Remark\ 3.} 
(1)\ \ 
We see that  $\sum_{l(Q) \leq 1}c(Q)\phi_{Q}$ is convergent in 
${\mathcal S}'$
for each  sequence $c \in a^s(e^{s'}_{pq})^{\sigma}_{x_{0}}$ or $c 
\in a^s(\tilde{e}^{s'}_{pq})^{\sigma}_{x_{0}}$ by Lemma 3.

(2)\ \ 
We notice  that $
D\equiv
\{ f= \sum_{l(Q)\leq 1}c(Q)\phi_{Q} : 
\ \ c \in a^s(e^{s'}_{pq})^{\sigma}_{x_{0}} \}$
 is independent of the choice of  $\phi_{0},\ \phi \in {\mathcal S}$ 
as in Remark 2.
Indeed, suppose $\{\phi_{0}^1, \phi^1\}$ and $\{\phi_{0}^2, \phi^2\}$ are Schwartz functions  as in Remark 2, and the spaces $D^1$ and $D^2$ 
are defined by using  $\{\phi_{0}^1, \phi^1\}$ and $\{\phi_{0}^2, \phi^2\}$ in the place of $\{\phi_{0},\phi\}$ respectively. We consider 
the $\varphi$-transform
$$
\phi^1_{P}=\sum_{l(R) \leq 1}l(R)^{-n}\langle \phi^1_{P}\ ,\ \ \varphi^2_{R} \rangle 
\phi^2_{R}.
$$
Then for $D^1 \ni 
f=\sum_{l(P)\leq 1}c(P)\phi^1_{P},\ c \in a^s(e^{s'}_{pq})^{\sigma}_{x_{0}}$, 
we have   
$$
f=\sum_{l(P) \leq 1}c(P)\phi^1_{P}=\sum_{l(R)\leq 1}Ac(R)\phi^2_{R}
$$
where $A=\{l(R)^{-n}\langle \phi^1_{P}\ ,\ \varphi^2_{R} \rangle \}_{RP}$.
From Lemma 1 and Lemma 2, we see that  for 
$c \in a^s(e^{s'}_{pq})^{\sigma}_{x_{0}}$,
$Ac \in a^s(e^{s'}_{pq})^{\sigma}_{x_{0}}$. 
This shows that
$D^1 \subset D^2$.
By the same argument, we see that  
$D^2 \subset D^1$. That is,  $D^1=D^2$. These imply that the space 
$D$ is independent of the choice of $\{\phi_{0}, \phi\}$.
In the same way $\tilde{D}=\{ f= \sum_{l(Q) \leq 1}c(Q)\phi_{Q} : 
\ \ (c(Q)) \in a^s(\tilde{e}^{s'}_{pq})^{\sigma}_{x_{0}} \}$ is
 independent of the choice of  $\{\phi_{0}, \phi\}$.

\bigskip

\noindent
{\it Proof of Theorem {\rm 1}.}\ \ (i)\ \ We may assume $\sigma \geq 0$ by Remark 1.  We put 
 $D \equiv \{ f= \sum_{l(Q) \leq 1}c(Q)\phi_{Q} : 
\ \ c \in a^s(e^{s'}_{pq})^{\sigma}_{x_{0}} \}$. 
In order to prove $D \subset A^s(E^{s'}_{pq})^{\sigma}_{x_{0}}$ we claim 
for a dyadic cube $P$, and for $f=\sum_{Q}c(Q)\phi_{Q} \in D$,
$$
c(E^{s'}_{pq})(P)\leq Cc(e^{s'}_{pq})(P)\ \ \ \ \ \ \ \ \ \ \ \ \ \ \ \ \ \ \ \ \ \ \ \ \ \ \ \ \ \ \ \ \ \ \ \ \ \ \ \ \ \ \ \ \ \ \ \ \ {\rm (a)}
$$
if $0 < p,q \leq \infty$. 
Let $(c(P)) \in a^s(e^{s'}_{pq})^{\sigma}_{x_{0}}$. 
Since $\mathcal{S}$ is closed under the convolution, 
we have, for $i \geq 0$,
\begin{eqnarray*}
\lefteqn{|\phi_{i} * f(x)| = 
|\sum_{l(P)\leq 1} c(P) \phi_{i}*\phi_{P} (x)|
}
\\
&=&|\sum_{j=(i-1)\vee 0}^{i+1}\sum_{l(P)=2^{-j}}c(P)\phi_{i}*\phi_{P}(x)|
\\
&\leq& 
C \sum_{j=(i-1)\vee 0}^{i+1}\sum_{l(P)=2^{-j}}|c(P)|(1+2^j|x-x_{P}|)^{-L} 
\end{eqnarray*}
for a sufficiently large number $L$.
Hence we have, 
using the maximal function $M_{t}f(x)$,\ $0< t\leq 1$,  
as in the proof of Lemma 1

\begin{eqnarray*}
\lefteqn{
\{\sum_{i\geq j\vee 0} (2^{is'}|\phi_{i}*f|)^{q}\}^{1/q} \leq 
C\{\sum_{i\geq j\vee 0} (2^{is'}\sum_{l(R)=2^{-i}}|\phi_{i}*f|\chi_{R})^q\}^
{1/q}
}
\\
&\leq& 
C\{ \sum_{i \geq j\vee 0}(2^{is'}
\sum_{l(R)=2^{-i}}
(\sum_{k=(i-1)\vee0}^{i+1}\sum_{l(R')=2^{-k}}|c(R')|
(1+2^{k}|x-x_{R'}|)^{-L})\chi_{R})^q\}^{1/q}
\\
&\leq&
C\{ \sum_{i \geq j\vee 0}(
\sum_{l(R)=2^{-i}}M_{t}(\sum_{k=(i-1)\vee0}^{i+1}\sum_{l(R')=2^{-k}}2^{is'}|c(R')|\chi_{R'})\chi_{R})^{q}\}^{1/q}
\end{eqnarray*}
if $0< t \leq 1$ and $L > n/t$. 
Taking $L^p(P)$-norm and using the  Fefferman-Stein inequality, 
we have  for a dyadic cube $P$ with $l(P)=2^{-j}$
\begin{eqnarray*}
\lefteqn{
c(F^{s'}_{pq})(P)=||\{\sum_{i\geq j\vee 0} (2^{is'}|\phi_{i}*f|)^{q}\}^{1/q}||_{L^p(P)} 
}
\\
&\leq&
C||\{ \sum_{i \geq j\vee 0}(
M_{t}(\sum_{k=(i-1)\vee0}^{i+1}\sum_{l(R')=2^{-k}}2^{is'}|c(R')|\chi_{R'}))^{q}\}^{1/q}||_{L^p(P)}
\\
&\leq&
C||\{ \sum_{i \geq j\vee 0}(
\sum_{k=(i-1)\vee0}^{i+1}\sum_{l(R')=2^{-k}}2^{is'}|c(R')|\chi_{R'}
)^{q}\}^{1/q}||_{L^p(P)}
\\
&\leq&
C||\{\sum_{i \geq j\vee 0}(
\sum_{l(R')=2^{-i}}2^{is'}|c(R')|
\chi_{R'})^{q}\}^{1/q}||_{L^p(P)}=Cc(f^{s'}_{pq})(P)
\end{eqnarray*}
if $0< t <\min(p,q)$ and
 $0< p < \infty$. For the $ p =\infty$ case , we obtain the same result. 
In the same way for the B-type case we have the same estimate 
$$
c(B^{s'}_{pq})(P)\leq Cc(b^{s'}_{pq})(P)\ \ 
{\rm if } \ \ 0< p \leq \infty.
$$ 
This implies 
$D \subset A^s(E^{s'}_{pq})^{\sigma}_{x_{0}}$. 

In order to complete the proof of Theorem 1 (i), we will show the inverse.
 We consider the 
$\varphi$-transform 
$f= \sum_{l(P)\leq 1}c(f)(P)\varphi_{P}$, 
$c(f)(P)= l(P)^{-n}\langle f\ ,\ \phi_{P} \rangle$ 
where $\phi_{P}$ and  $\varphi_{P}$ as in Remark 2. 
It suffices to show  that $c(f)(P) \in a^s(e^{s'}_{pq})^{\sigma}_{x_{0}}$ 
for $f \in A^s(E^{s'}_{pq})^{\sigma}_{x_{0}}$. 
More precisely, we claim that for a dyadic cube $P$ with $l(P)=2^{-j}$,

$$
c(f)(e^{s'}_{pq})(P) \leq Cc(E^{s'}_{pq})(P)\ \ \ \ \ \ \ \ \ \ \ \ \ \ \ \ \ \ \ \ \ \ \ \ \ \ \ \ \ \ \ \ \ \ \ \ \ \ \ \ \ \ \ \ \ \ {\rm (b)}
$$ 
where $c(f)(e^{s'}_{pq})(P)$ is a sequence defined by  replacing 
 the sequence $c(P)$ by the sequence $c(f)(P)$ 
in the definition of $c(e^{s'}_{pq})(P)$. For $f \in \mathcal{S}'$ and a dyadic cube $P$ with $l(P)=2^{-j}$,
we define the sequence $\sup(f)(P)$ by setting 
$$
\sup(f)(P)= 
\sup_{P \ni y}|\phi_{j}*f(y)|.
$$  For $\gamma \in \mathbb{N}_{0}$ the sequences $\inf_{\gamma}(f)(P)$, 
$t_{\gamma}(P)$ are defined previously and  for a sequence $c(P)$, we also define a sequence $c^*(P)$ previously (See  Lemma 4).
We have, from the fact in [9, Lemma A.4] that
$
\sup(f)^*(P) \sim \inf_{\gamma}(f)^*(P)
$
for $\gamma$ sufficiently large. 

Thus, we have

\bigskip

\noindent
$
|c(f)(P)|= l(P)^{-n}|\langle f\ ,\ \phi_{P} \rangle| 
=|\phi_{j} * f (x_{P})| 
\leq
 \sup(f)(P) 
\leq \sup(f)^*(P) 
\sim \inf_{\gamma}(f)^*(P)
$

\bigskip
\noindent
for $\gamma$ sufficiently large.  
 Therefore, from Lemma 4 (i) and (ii) we have 

\begin{eqnarray*}
\lefteqn{|c(f)(f^{s'}_{pq})(P)| \leq C {\rm inf}_{\gamma}(f)^*(f^{s'}_{pq})(P) \leq C {\rm inf}_{\gamma}(f)(f^{s'}_{pq})(P)
}
\\
&\leq&
C||\{ \sum_{i\geq j\vee 0}(\sum_{l(R)=2^{-i}}2^{is'}{\rm inf}_{\gamma}(f)(R)
\chi_{R})^{q} \}^{1/q}||_{L^p(P)} 
\\
&\leq&
C||\{ \sum_{i\geq j\vee 0}(2^{is'}2^{\gamma L}
\sum_{l(R')=2^{-(\gamma+i)}}t_{\gamma}^*(R')\chi_{R'})^{q} \}^{1/q}||_{L^p(P)} 
\\
&\leq&
C2^{\gamma L}||\{ \sum_{i\geq (j\vee 0)+\gamma}(2^{is'}2^{-\gamma s'}\sum_{l(R')=2^{-i}}t_{\gamma}(R')\chi_{R'})^{q} \}^{1/q}||_{L^p(P)} 
\\
&\leq&
C2^{\gamma( L- s')}||\{ \sum_{i\geq (j\vee 0)+\gamma}(2^{is'}\sum_{l(R')=2^{-i}}
|\phi_{i-\gamma}*f(y)|\chi_{R'})^{q} \}^{1/q}||_{L^p(P)} 
\\
&\leq&
C2^{\gamma( L- s')}||\{ \sum_{i\geq j\vee 0}(2^{is'}2^{s'\gamma}
\sum_{l(R')=2^{-(i+\gamma)}}
|\phi_{i}*f(y)|\chi_{R'})^{q} \}^{1/q}||_{L^p(P)} 
\\
&\leq&
C2^{\gamma L}||\{ \sum_{i\geq j\vee 0}(2^{is'}
|\phi_{i}*f(y)|)^{q} \}^{1/q}||_{L^p(P)}=Cc(F^{s'}_{pq})(P) 
\end{eqnarray*}
if $0< p < \infty$. For $ p =\infty$ , we obtain the same result. 
For the B-type case we can prove the same result by the same argument as above,
$$
c(f)(b^{s'}_{pq})(P) \leq C c(B^{s'}_{pq})(P)
$$ 
if $ 0 < p \leq \infty$. Thus, we obtain 
$$
c(f)(e^{s'}_{pq})(P) \leq C c(E^{s'}_{pq})(P).
$$
By Remark 3 (2) this implies that,
$A^s(E^{s'}_{pq})^{\sigma}_{x_{0}} \subset D$, $0< p \leq \infty$. 
Hence, we obtain $A^s(E^{s'}_{pq})^{\sigma}_{x_{0}} = D$.

(ii)\ \ We can prove (ii) in the same way as  (i).
\qed

\bigskip

We have the following properties from Theorem 1.

\bigskip

\noindent
{\bf Proposition 1.}\ \ {\it 
 Suppose that 
\ $s, \ s,' \ \sigma \in {\mathbb R}$ and $x_{0} \in \mathbb{R}^n$.

{\rm (i)}\ \ When $\sigma < 0$, we have 
\ \ 
$A^s(E^{s'}_{pq})^{\sigma}_{x_{0}}= \{ 0 \}$,\ \ for\ \ $0 < p, q \leq \infty$,

{\rm (ii)}\ \ When $\sigma+s < 0$, 
we have 
$A^s(B^{s'}_{pq})^{\sigma}_{x_{0}}= \{ 0 \}$,\ \  for\ \ $0 < p,q \leq \infty$, \ \ 
and $A^s(F^{s'}_{pq})^{\sigma}_{x_{0}}= \{ 0 \}$,\ \  for\ \ $0 < p < \infty$,\  $0<q \leq \infty$,
 
(iii)\ \  When $s < 0$, 
we have 
$A^s(\tilde{B}^{s'}_{pq})^{\sigma}_{x_{0}}= \{ 0 \}$,\ \  for\ \ $0 < p,q \leq \infty$,
\ \ 
and 
$A^s(\tilde{F}^{s'}_{pq})^{\sigma}_{x_{0}}= \{ 0 \}$,\ \  for\ \ $0 < p < \infty,\  0<q \leq \infty$. 
}
\bigskip

\noindent
{\it Proof.}\ \ These properties are shown easily.
\qed

\bigskip

\noindent
{\bf Proposition 2.}\ \ {\it 
 Suppose that 
\ $s, \ s,' \ \sigma \in {\mathbb R}$ and $x_{0} \in \mathbb{R}^n$.

\begin{enumerate}[{\rm (i)}]
\item \ 
When $s \leq 0$, we have 

$A^s(B^{s'}_{pq})^{\sigma}_{x_{0}}= (B^{s'}_{pq})^{s+\sigma}_{x_{0}}$\ 
for\ \ $0 < p,q \leq \infty$,
and
$A^s(F^{s'}_{pq})^{\sigma}_{x_{0}}= (F^{s'}_{pq})^{s+\sigma}_{x_{0}}$\ 
for\ \ $0 < p < \infty,\ 0<q \leq \infty$, 

In particular,
  when $\sigma \geq 0$ and $\sigma+s = 0$, 
we have 

$A^s(B^{s'}_{pq})^{\sigma}_{x_{0}}= B^{s'}_{pq}(\mathbb{R}^n)$ \ for\ \ $0 < p,q \leq \infty$, 
and
$A^s(F^{s'}_{pq})^{\sigma}_{x_{0}}= F^{s'}_{pq}(\mathbb{R}^n)$ \ 
for\ \ $0 < p <\infty,\ 0<q \leq \infty$. 

\item \ When $\sigma \geq 0$, we have 

$A^s(E^{s'+\sigma}_{pq}) \subset A^s(\tilde{E}^{s'}_{pq})^{\sigma}_{x_{0}} \subset  A^s(E^{s'}_{pq})^{\sigma}_{x_{0}}$,

and when $\sigma < 0$, 
we have 

$A^s(\tilde{E}^{s'}_{pq})^{\sigma}_{x_{0}} \subset A^s(E^{s'+\sigma}_{pq})$. 

\item \ \ If $\sigma \geq 0$, then we have 

$(E^{s'}_{\infty \infty})^{\sigma}_{x_{0}}=
(\tilde{E}^{s'}_{\infty\infty})^{\sigma}_{x_{0}}$.

\end{enumerate}
}

\bigskip

\noindent
{\it Proof.}\ \ The property (i) can be proved from the fact that 
$$
Cl(Q)^{-\sigma}
\sup_{\mathcal{D} \ni P \subset 3Q}l(P)^{-s}c(e^{s'}_{pq})(P)
\geq l(Q)^{-(\sigma+s)}c(e^{s'}_{pq})(Q)
$$
and
$$
l(Q)^{-\sigma}\sup_{\mathcal{D} \ni P \subset 3Q}
l(P)^{-s}c(e^{s'}_{pq})(P)
\leq  Cl(Q)^{-(\sigma+s)}\sup_{\mathcal{D} \ni P \subset 3Q}c(e^{s'}_{pq})(P),
$$
if $s \leq 0$.

We obtain the property (ii) from the fact that  
$$
c(\tilde{e}^{s'}_{pq})^{\sigma}_{x_{0}}(P) \leq Cc(e^{\sigma+s'}_{pq})(P), 
$$
and 
$$
l(Q)^{-\sigma}l(P)^{-s}c(e^{s'}_{pq})(P) \leq Cl(P)^{-s}c(\tilde{e}^{s'}_{pq})^{\sigma}_{x_{0}}(P)
$$
since $l(Q)^{-\sigma} \leq C(l(P)+|x_{0}-x_{P}|)^{-\sigma}$ for $P \subset 3Q$ if $\sigma \geq 0$.
The last half of property (ii) can be proved  since 
$$
c(\tilde{e}^{s'}_{pq})^{\sigma}_{x_{0}}(P)
 \geq c(e^{s'+\sigma}_{pq})(P) 
$$
if $\sigma < 0$.
To prove the property (iii), it suffices to see from property (ii),

$$
(e^{s'}_{\infty \infty})^{\sigma}_{x_{0}} \subset
  (\tilde{e}^{s'}_{\infty\infty})^{\sigma}_{x_{0}}.
$$  

\noindent 
 We consider  any dyadic cube $R$ with $l(R)=2^{-i}$ and 
 dyadic cubes $Q_{l}$ with $x_{0} \in Q_{l}$ and $l(Q_{l})=2^{-l},\ i \geq l$ such that 
 $Q_{i} \subset \cdots  \subset Q_{l} \subset Q_{l-1} \subset \cdots $ and 
 $\cup_{i \geq l}Q_{l}=\mathbb{R}^n$. We set 
 $Q_{l}^0 \equiv 3Q_{l} \setminus 3Q_{l+1},\ i >l$ 
 and $Q_{i}^0 \equiv 3Q_{i}$. 
 We divide the proof into two cases: 
 
 Case (a): $R \subset Q_{l}^0,\ i > l$ case. Then we have 
 $2^{-i}+|x_{0}-x_{R}| \geq C 2^{-l},$
 
 Case (b):  $R \subset Q_{i}^0$. Then we have 
 $2^{-i}+ |x_{0} -x_{R}| \geq 2^{-i}$. 
 
\noindent  
In the case (a) we have
\begin{eqnarray*}\lefteqn{
2^{is'}|c(R)|(2^{-i}+|x_{0}-x_{R}|)^{-\sigma} \leq C 2^{is'}2^{l\sigma}|c(R)| 
}
\\
&&
\leq C \sup_{x_{0}\in Q}2^{l\sigma}\sup_{R \subset 3Q}2^{is'}|c(R)| < \infty.
\end{eqnarray*}
In the case (b) we have 
\begin{eqnarray*}\lefteqn{
 2^{is'}|c(R)|(2^{-i}+|x_{0}-x_{R}|)^{-\sigma}  \leq C 2^{is'}2^{i\sigma}|c(R)| }
\\
&&
\leq C
 \sup_{x_{0} \in Q}2^{l\sigma}\sup_{R \subset 3Q}2^{is'}|c(R)|< \infty.
\end{eqnarray*}
The proof is complete.\qed

\bigskip

\noindent
{\bf Proposition 3.}\ \ {\it 
 Suppose that 
\ $s, \ s,' \ \sigma \in {\mathbb R}$, and $x_{0} \in \mathbb{R}^n$.

When  $0 < q_{1} \leq q_{2} \leq \infty$, $0 <p \leq \infty$,
we have 

$$
A^{s}(B^{s'}_{pq_{1}})^{\sigma} _{x_{0}}
\subset 
A^s(B^{s'}_{pq_{2}})^{\sigma}_{x_{0}},\ 
\ A^{s}(\tilde{B}^{s'}_{pq_{1}})^{\sigma} _{x_{0}}
\subset 
A^s(\tilde{B}^{s'}_{pq_{2}})^{\sigma}_{x_{0}}
$$,

 and  
when  $0 < q_{1} \leq q_{2} \leq \infty$, $0 <p < \infty$,
 we have 

$$
A^{s}(F^{s'}_{pq_{1}})^{\sigma} _{x_{0}}
\subset 
A^s(F^{s'}_{pq_{2}})^{\sigma}_{x_{0}},\  
\ A^{s}(\tilde{F}^{s'}_{pq_{1}})^{\sigma} _{x_{0}}
\subset 
A^s(\tilde{F}^{s'}_{pq_{2}})^{\sigma}_{x_{0}}
$$.

}
\bigskip

\noindent
{\it Proof}.\ \ These  inclusions are  corollaries of the monotonicity 
of the $l^{p}$-norm.\qed

\bigskip

\noindent
{\bf Proposition 4.}\ \ {\it 
 Suppose that 
\ $s, \ s,' \ \sigma \in {\mathbb R}$,\ $0 <\epsilon$ and $x_{0} \in \mathbb{R}^n$. We have

\noindent
\begin{enumerate}[{\rm (i)}]

\item\ 
 $A^{s}(B^{s'+\epsilon}_{pq_{1}})^{\sigma-\epsilon}_{x_{0}} 
\subset 
A^{s}(B^{s'}_{pq_{2}})^{\sigma}_{x_{0}}$ for\ \  $0<p\leq \infty,\ 0< q_{1}, \ q_{2} \leq \infty$, and

$A^{s}(F^{s'+\epsilon}_{pq_{1}})^{\sigma-\epsilon}_{x_{0}} 
\subset 
A^{s}(F^{s'}_{pq_{2}})^{\sigma}_{x_{0}}$ for\ \  $0<p< \infty,\ 0< q_{1}, \ q_{2} \leq \infty$, and

\item\ 
$A^{s+\epsilon}(E^{s'}_{pq})^{\sigma-\epsilon}_{x_{0}} 
\subset 
A^{s}(E^{s'}_{pq})^{\sigma}_{x_{0}}$ for\ \  $0< p, \ q \leq \infty$, and

\item\ \ \ 
 $A^{s-\epsilon}(B^{s'+\epsilon}_{pq_{1}})^{\sigma}_{x_{0}} 
\subset 
A^{s}(B^{s'}_{pq_{2}})^{\sigma}_{x_{0}}$, and 
$A^{s-\epsilon}(\tilde{B}^{s'+\epsilon}_{pq_{1}})^{\sigma}_{x_{0}} 
\subset 
A^{s}(\tilde{B}^{s'}_{pq_{2}})^{\sigma}_{x_{0}}$ 
for\ \  $0< p,\ q_{1}, \ q_{2} \leq \infty$, and

 $A^{s-\epsilon}(F^{s'+\epsilon}_{pq_{1}})^{\sigma}_{x_{0}} 
\subset 
A^{s}(F^{s'}_{pq_{2}})^{\sigma}_{x_{0}}$, 
 and 
$A^{s-\epsilon}(\tilde{F}^{s'+\epsilon}_{pq_{1}})^{\sigma}_{x_{0}} 
\subset 
A^{s}(\tilde{F}^{s'}_{pq_{2}})^{\sigma}_{x_{0}}$ 
for\ \  $0< p < \infty,\ 0< q_{1}, \ q_{2} \leq \infty$.

\end{enumerate}

}

\bigskip

\noindent
{\it Proof}.\ \ (ii) is obvious.
 (i) and (iii) are corollaries of  H$\ddot{\rm o}$lder's inequality and the monotonicity of the $l^{p}$-norm.
\qed

\bigskip

\noindent
{\bf Proposition 5.}\ \ {\it 
 Suppose that 
\ $s, \ s,' \ \sigma \in {\mathbb R}$ and $x_{0} \in \mathbb{R}^n$.

\begin{enumerate}[{\rm (i)}]
\item\ If 
$ 0 < p_{2} \leq p_{1} \leq \infty$ and $0 < q \leq \infty$, then 

$A^{s+\frac{n}{p_{1}}}(B^{s'}_{p_{1}q})^{\sigma}_{x_{0}}
\subset 
A^{s+\frac{n}{p_{2}}}(B^{s'}_{p_{2}q})^{\sigma}_{x_{0}}$, 
$A^{s+\frac{n}{p_{1}}}(\tilde{B}^{s'}_{p_{1}q})^{\sigma}_{x_{0}}
\subset 
A^{s+\frac{n}{p_{2}}}(\tilde{B}^{s'}_{p_{2}q})^{\sigma}_{x_{0}}$,

 and,
if 
$ 0 < p_{2} \leq p_{1} < \infty$ and $0 < q \leq \infty$, then 

$A^{s+\frac{n}{p_{1}}}(F^{s'}_{p_{1}q})^{\sigma}_{x_{0}}
\subset 
A^{s+\frac{n}{p_{2}}}(F^{s'}_{p_{2}q})^{\sigma}_{x_{0}}$, 
$A^{s+\frac{n}{p_{1}}}(\tilde{F}^{s'}_{p_{1}q})^{\sigma}_{x_{0}}
\subset 
A^{s+\frac{n}{p_{2}}}(\tilde{F}^{s'}_{p_{2}q})^{\sigma}_{x_{0}}$. 

\item\ 
If $ 0<q \leq \infty,\ 0 < p \leq \infty,\ \frac{n}{p} < s$, 
then

$A^{s}(E^{s'}_{pq})^{\sigma}_{x_{0}}=
(E^{s+s'-\frac{n}{p}}_{\infty\infty})^{\sigma}_{x_{0}}$ and
$A^{s}(\tilde{E}^{s'}_{pq})^{\sigma}_{x_{0}}=
(\tilde{E}^{s+s'-\frac{n}{p}}_{\infty\infty})^{\sigma}_{x_{0}}$.

In particular, if $0 \leq \sigma ,\  0<q \leq \infty,\ 0 < p \leq \infty,
\ \frac{n}{p} < s$, 
then

$A^{s}(E^{s'}_{pq})^{\sigma}_{x_{0}}=
A^{s}(\tilde{E}^{s'}_{pq})^{\sigma}_{x_{0}}$.

\item\ If \ \ 
$ 0<  p_{1},\ p_{2},\ q \leq \infty$, then

 $A^{\frac{n}{p_{1}}}(E^{s'}_{p_{1}\infty})^{\sigma}_{x_{0}}
= 
A^{\frac{n}{p_{2}}}(E^{s'}_{p_{2}\infty})^{\sigma}_{x_{0}}
=(E^{s'}_{\infty \infty})^{\sigma}_{x_{0}}$,

  $A^{\frac{n}{p_{1}}}(\tilde{E}^{s'}_{p_{1}\infty})^{\sigma}_{x_{0}}
= 
A^{\frac{n}{p_{2}}}(\tilde{E}^{s'}_{p_{2}\infty})
^{\sigma}_{x_{0}}
=(\tilde{E}^{s'}_{\infty \infty})^{\sigma}_{x_{0}}$,

 $A^{\frac{n}{p_{1}}}(F^{s'}_{p_{1}q})^{\sigma}_{x_{0}}
= 
A^{\frac{n}{p_{2}}}(F^{s'}_{p_{2}q})
^{\sigma}_{x_{0}}$,
  $A^{\frac{n}{p_{1}}}(\tilde{F}^{s'}_{p_{1}q})^{\sigma}_{x_{0}}
= 
A^{\frac{n}{p_{2}}}(\tilde{F}^{s'}_{p_{2}q})
^{\sigma}_{x_{0}}$.
\end{enumerate}
}
\bigskip

\noindent
{\it Proof}.\ \ The properties (i) are  corollaries of 
H$\ddot{\rm {o}}$lder's inequality. We will prove the properties  (ii).
We see that 
$$
a^{s+\frac{n}{p}}(e^{s'}_{pq})^{\sigma}_{x_{0}} \subset
(e^{s'+s}_{\infty\infty})^{\sigma}_{x_{0}},
$$
since
$$
l(P)^{-(s+\frac{n}{p})}c(e_{pq}^{s'})(P)\geq l(P)^{-(s'+s)}|c(P)|.
$$
 Hence in order to prove (ii), it suffices to prove  
$$
(e^{s'+s}_{\infty\infty})^{\sigma}_{x_{0}} \subset 
a^{s+\frac{n}{p}}(e^{s'}_{pq})^{\sigma}_{x_{0}}.
$$

\noindent
Since 
$$
c(\dot{e}^{s'}_{pq})(P) \leq C (e^{s'+s}_{\infty\infty})(P)\times 
l(P)^{s+\frac{n}{p}}
$$ 

\noindent
if $s > 0$ and $0 < q < \infty$, 
we get the desired result. Similarly, for the other case, we can prove.
  
The first part of  properties (iii) is  obtained in the same way in the proof of (ii) and the last part  is just [ 10: Corollary 5.7 ].  \qed

\bigskip

\noindent
{\bf Proposition 6.}\ (Embedding)\ \ {\it 
 Let $s,\ s',\ \sigma  \in {\Bbb R},\  0 < p,\ q \leq \infty\ and 
\ x_{0} \in \mathbb{R}^n$. We have 
\begin{enumerate}[{\rm (i)}]
\item\ \ 
$A^s(E^{s_{1}'}_{p \xi})^
{\sigma}_{x_{0}} \subset  A^s(E^{s_{2}'}_{p \eta})^{\sigma}_{x_{0}}$,
$A^s(\tilde{E}^{s_{1}'}_{p \xi})^
{\sigma}_{x_{0}} \subset  A^s(\tilde{E}^{s_{2}'}_{p \eta})^{\sigma}_{x_{0}}$,
for $s_{1}'> s_{2}'$ and $0< \xi, \eta \leq \infty$,

\item \ \ 
$A^s(B^{s_{1}'}_{p_{1}q})^{\sigma}_{x_{0}}
 \subset  A^s(B^{s_{2}'}_{p_{2} q})^{\sigma}_{x_{0}}$,\ \  
$A^s(\tilde{B}^{s_{1}'}_{p_{1}q})^{\sigma}_{x_{0}} \subset 
 A^s(\tilde{B}^{s_{2}'}_{p_{2}q})^{\sigma}_{x_{0}}$, 
for $s_{1}'-s_{2}' =
n(\frac{1}{p_{1}}-\frac{1}{p_{2}})$ 
and 
$ 0< p_{1} \leq p_{2} \leq \infty$,

$A^s(F^{s_{1}'}_{p_{1}\xi})^{\sigma}_{x_{0}}
 \subset  A^s(F^{s_{2}'}_{p_{2} \eta})^{\sigma}_{x_{0}}$,\ \  
$A^s(\tilde{F}^{s_{1}'}_{p_{1}\xi})^{\sigma}_{x_{0}} \subset 
 A^s(\tilde{F}^{s_{2}'}_{p_{2}\eta})^{\sigma}_{x_{0}}$, 
for $s_{1}'-s_{2}' =
n(\frac{1}{p_{1}}-\frac{1}{p_{2}})$ 
and 
$ 0< p_{1} < p_{2} < \infty$, $0 < \xi, \eta \leq \infty$,

\item \ \ 
$A^s(B^{s'}_{pq})^{\sigma}_{x_{0}}
 \subset  A^s(F^{s'}_{p q})^{\sigma}_{x_{0}}$,\ \  
$A^s(\tilde{B}^{s'}_{pq})^{\sigma}_{x_{0}} \subset 
 A^s(\tilde{F}^{s'}_{pq})^{\sigma}_{x_{0}}$, 
for $0< q \leq p \leq \infty$,

$A^s(F^{s'}_{pq})^{\sigma}_{x_{0}}
 \subset  A^s(B^{s'}_{p q})^{\sigma}_{x_{0}}$,\ \  
$A^s(\tilde{F}^{s'}_{pq})^{\sigma}_{x_{0}} \subset 
 A^s(\tilde{B}^{s'}_{pq})^{\sigma}_{x_{0}}$, 
for $0< p \leq q \leq \infty$.
\end{enumerate}
}

\bigskip

\noindent
{\it Proof}.\ \ 
 The  embedding properties (i) and  the first embedding of (ii) are
 corollaries of  H$\ddot{\rm o}$lder's inequality and the
monotonicity property of the $l^{p}$-norm. For the second embedding of (ii), see [37; Proposition 2.5] (cf. [38; Theorem 2.7.1]).
 (iii) is a  corollary of Minkowski's inequality
 (cf. Triebel[ 38: 2.3.2 Proposition 2 ] ).

\qed

\bigskip

\noindent
{\bf Remark 4.}\ \ 
 Let $0 < p,q \leq \infty$,\ $s, \sigma \in \mathbb{R}$,\ $x_{0} \in \mathbb{R}^n$ and $s'> n(\frac{1}{p}-1)_{+}$. If $f \in A^s(E^{s'}_{pq})^{\sigma}_{x_{0}}$, then $f$ is locally integrable (and  locally $L^p$ integrable). Indeed, 
we consider the Littlewood-Paley decomposition  
$$
f=\sum_{i \geq 0}f*\phi_{i}.
$$
It suffices to show that
$\sum_{i \geq 0}f*\phi_{i}$ is locally integrable and locally $L^p$ integrable.
We may consider any dyadic cube $P$ with $l(P) \geq  1$.  Then we have   
if $1 \leq p < \infty$,
\begin{eqnarray*}
\lefteqn{||\sum_{i \geq  0}f*\phi_{i}||_{L^1(P)} 
\leq C||\sum_{i \geq  0}f*\phi_{i}||_{L^p(P)} 
}
\\
&\leq& C||\{\sum_{i \geq  0}(2^{is'}|f*\phi_{i}|)^q\}^{1/q}||_{L^p(P)}
\leq C c(F^{s'}_{pq})(P)   < \infty
\end{eqnarray*}
by using H${\rm \ddot{o}}$lder inequality if $1 \leq q \leq \infty$ and the
monotonicity property of the $l^{p}$-norm if $0 < q \leq 1$.  In the same way we have 
\begin{eqnarray*}
\lefteqn{||\sum_{i \geq  0}f*\phi_{i}||_{L^1(P)} 
 \leq C||\sum_{i \geq  0}f*\phi_{i}||_{L^p(P)} \leq C\sum_{i \geq 0}||f*\phi_{i}||_{L^p(P)}
}
\\
&\leq& C\{\sum_{i \geq  0}(2^{is'}||f*\phi_{i}||_{L^p(P)})^q\}^{1/q}
\leq C c(B^{s'}_{pq})(P)   < \infty. 
\end{eqnarray*}
If $0 < p \leq 1$, in the same way we have 
\begin{eqnarray*}
\lefteqn{||\sum_{i \geq  0}f*\phi_{i}||_{L^p(P)} \leq C||\sum_{i \geq   0}f*\phi_{i}||_{L^1(P)} 
}
\\
&\leq& C\{\sum_{i \geq  0}(2^{i(s'-n(\frac{1}{p}-1))}||f*\phi_{i}||_{L^1(P)})^q\}^{1/q}
\\
&=&C c(B^{s'-n(\frac{1}{p}-1)}_{1q})(P)\leq C c(B^{s'}_{pq})(P) < \infty, 
\end{eqnarray*}
where we use Proposition 6  in the last inequality.
Similarly,  by using the fact that 
$c(B^{s'}_{p\ p\vee q})(P) \leq c(F^{s'}_{pq})(P)$ 
 we have the same estimate for the F-type case  if $0 < p \leq 1$. 
Therefore, we obtain  the desired result for $f \in A^s(E^{s'}_{pq})^{\sigma}_{x_{0}}$.  But
we note that it holds for $0 < p < \infty$ in the F-type case  and  for 
$0 < p \leq  \infty$ in the B-type case.  
We note that it holds an analogous result 
for $f \in A^s(\tilde{E}^{s'}_{pq})^{\sigma}_{x_{0}}$ with the weight $w_{i}
=(2^{-i}+|x_{0}-x|)^{-\sigma}$.
\bigskip

We recall the definitions  of smooth atoms and molecules. 

\bigskip

\noindent
{\bf Definition 6.} \ \ 
 Let $r_{1},\ r_{2} \in \mathbb{N}_{0}, L >n $. A family of functions $m = (m_{Q})$ indexed by dyadic cubes $Q$ 
with 
$l(Q) \leq 1$
is called a family of $(r_{1}, r_{2}, L)$-- smooth molecules  if 

\bigskip

$(3.1)$ \ \ 
 $
|m_{Q}(x)| \leq C (1+l(Q)^{-1}|x-x_{Q}|)^{-\max(L, L_{2})}
$ for some $L_{2} > n+r_{2}$ when $l(Q) < 1$,

$(3.2)$ \ \ 
$|\partial^{\gamma}m_{Q}(x)| \leq Cl(Q)^{-|\gamma|}(1+l(Q)^{-1}|x-x_{Q}|)
^{-L}
$
 for $0 <|\gamma| \leq r_{1}$, when $l(Q) < 1$ and

$(3.3)$\ \  $\int_{{\mathbb R}^n}  x^{\gamma}m_{Q}(x) dx =0$ for $|\gamma| < r_{2}$ 
when $l(Q) < 1$,

\noindent
where $(3.2)$ is void when $r_{1}=0$, and (3.3) is void when $r_{2}=0$,

$(3.4)$\ \  $|\partial^{\gamma}m_{Q}(x)| \leq Cl(Q)^{-|\gamma|}
(1+l(Q)^{-1}|x-x_{Q}|)^{-L},\ \ |\gamma| \leq r_{1}$ when $l(Q)=1$,

$(3.5)$\ \  we do not assume the vanishing moment condition (3.3) when $l(Q) =1$. 

\noindent
A family of functions $a = (a_{Q})$ indexed by dyadic cubes $Q$ with $l(Q) \leq 1$ 
is  called a family of  $(r_{1}, r_{2})$-- smooth atoms 
if 

\bigskip

$(3.6)$ \ \ 
supp $a_{Q} \subset 3Q$ for each dyadic cube $Q$ when $l(Q) \leq 1$,

$(3.7)$ \ \ 
 $
|\partial^{\gamma}a_{Q}(x)| \leq Cl(Q)^{-|\gamma|}
$
 for $|\gamma| \leq r_{1}$ when $l(Q) \leq 1$, and

 $(3.8)$\ \  $\int_{{\mathbb R}^n}  x^{\gamma}a_{Q}(x) dx =0$ for $|\gamma| < r_{2}$ 
when $l(Q) < 1$,

\noindent
where $(3.8)$ is void when $r_{2}=0$,

$(3.9)$\ \ we do not assume the vanishing moment condition (3.8)
when $l(Q)=1$.

\bigskip

\noindent
{\bf Theorem 2.}\ \ {\it 
Let $s, \ s',\  \sigma \in {\Bbb R}$,\ $0 < p,\ q \leq \infty$\ and $x_{0} \in \mathbb{R}^n$.
Let $r_{1}, r_{2} \in \mathbb{N}_{0}$, \  $J$  as in Lemma {\rm 1}  
and $L > n$.

{\rm (i)}\ We assume that  $r_{1}$, $r_{2}$ and $L$ satisfy the following conditions
 
\bigskip

 {\rm (4.1)}\ \ $r_{1} > \max (s',\ \sigma+s + s'- \frac{n}{p})$, 

 {\rm (4.2)}\ \ $r_{2} > J-n-s'$,

  {\rm (4.3)}\  \ $L > J$. 

Then we have

\begin{eqnarray*}
\lefteqn{
A^s(E^{s'}_{pq})^{\sigma}_{x_{0}}
=
 \{ f =\sum_{l(Q) \leq 1}c(Q)m_{Q} :
}
\\ 
\lefteqn{  
(r_{1},r_{2},L){\rm -\  smooth\ molecules}\ (m_{Q}),
\ \ (c(Q)) \in a^s(e^{s'}_{pq})^{\sigma}_{x_{0}} \}
}
\\
&=&
\{f= \sum_{l(Q)\leq 1}c(Q)a_{Q} :
\\
&&
   (r_{1},r_{2}){\rm - \ smooth\  atoms}\ (a_{Q}),
\ \ (c(Q)) \in a^s(e^{s'}_{pq})^{\sigma}_{x_{0}} \}.
\end{eqnarray*}

{\rm (ii)}\ 
 We assume that  $r_{1}$, $r_{2}$ and $L$ satisfy 
 
 \bigskip

 {\rm $(4.1)'$}\ \ 
$r_{1} > \max (s'+(\sigma\vee 0),\ (\sigma\vee 0)+s + s'- \frac{n}{p})$,

 {\rm $(4.2)'$}\ \ $r_{2} > J-n-s'-(\sigma\wedge 0)$,

 {\rm $(4.3)'$}\ \ $L > J+\sigma$

\noindent
Then we have

\begin{eqnarray*}
\lefteqn{
A^s(\tilde{E}^{s'}_{pq})^{\sigma}_{x_{0}}
=
\{ f =\sum_{l(Q) \leq 1}c(Q)m_{Q} :
}
\\
\lefteqn{
(r_{1},r_{2},L)
{\rm - \ smooth\ molecules}\ (m_{Q}),
\ \ (c(Q)) \in a^s(\tilde{e}^{s'}_{pq})^{\sigma}_{x_{0}} \}
}
\\
&=&
\{f= \sum_{l(Q)\leq 1}c(Q)a_{Q} :
\\&&
  (r_{1},r_{2})
{\rm - \ smooth\  atoms}\ (a_{Q}),
\ \ (c(Q)) \in a^s(\tilde{e}^{s'}_{pq})^{\sigma}_{x_{0}} \}.
\end{eqnarray*}

}

 \bigskip

\noindent
{\bf Remark 5.}\ \    
From  Lemma 3, 
we remark that  
$f =\sum_{l(Q) \leq 1}c(Q)m_{Q}$ and 
$f =\sum_{l(Q) \leq 1}c(Q)a_{Q}$ are convergent in 
 ${\mathcal S}'$ 
for each $c \in a^s(e^{s'}_{pq})^{\sigma}_{x_{0}}$ or 
$a^s(\tilde{e}^{s'}_{pq})^{\sigma}_{x_{0}}$.

\bigskip

\noindent
{\it Proof of Theorem {\rm 2}.}\ \  (i)\ \ We may assume $\sigma \geq 0$ by Remark 1. We put 

\noindent
$A\equiv \{f= \sum_{l(Q) \leq 1}c(Q)a_{Q} 
: \ \ (r_{1},r_{2}){\rm - smooth\  atoms}\ (a_{Q}),
 \ (c(Q)) \in a^s(e^{s'}_{pq})^{\sigma}_{x_{0}} \}$, \\
$M \equiv \{ f =\sum_{l(Q) \leq 1}c(Q)m_{Q} : \ \ 
(r_{1},r_{2},L){\rm - smooth\ molecules}\ (m_{Q}),\\ 
 (c(Q)) \in a^s(e^{s'}_{pq})^{\sigma}_{x_{0}} \}$.

Since an $(r_{1}, r_{2})$-- atom is an $(r_{1}, r_{2}, L)$-- molecule, it is easy to see that   $A \subset M$.
Let $M \ni f = \sum_{l(Q) \leq 1} c(Q) m_{Q}$ and  we consider the 
$\varphi$-transform 
$$
m_{Q}=\sum_{l(P)\leq 1} l(P)^{-n}\langle m_{Q}\ ,\ \varphi_{P} \rangle \phi_{P},
$$
where $\phi_{P}$ and $\varphi_{P}$ as in Remark 2.
Then we have

$$
f=\sum_{l(Q) \leq 1} c(Q) m_{Q} = \sum_{l(P)\leq 1} (Ac)(P) \phi_{P},
$$
where $A=\{ l(P)^{-n}\langle m_{Q}\ ,\ \varphi_{P} \rangle \}_{PQ}$.
Lemma 1 and Lemma 2 yield that $A$ is $(r_{1},\ r_{2}+n,\ L)$-- almost diagonal and 
$Ac \in a^s(e^{s'}_{pq})^{\sigma}_{x_{0}}$
for $c \in a^s(e^{s'}_{pq})^{\sigma}_{x_{0}}$.
Hence, if we put $D \equiv \{ f= \sum_{l(Q) \leq 1}c(Q)\phi_{Q} : 
\ \ c \in a^s(e^{s'}_{pq})^{\sigma}_{x_{0}} \}$, then  we see that 
$M \subset  D$. 
From Theorem 1 we see  $D = A^s(E^{s'}_{pq})^{\sigma}_{x_{0}}$. 
Hence, we obtain $A \subset M \subset 
A^s(E^{s'}_{pq})^{\sigma}_{x_{0}}$.

Using the argument similar to the proof of  
[10: Theorem 4.1] (cf. [4: Theorem 5.9] or [5: Theorem 5.8]), 
for $D \ni f=\sum_{l(Q) \leq 1}c(Q)\phi_{Q},\ \  c \in 
a^s(e^{s'}_{pq})^{\sigma}_{x_{0}}$,  
we see that there exist a family of $(r_{1}, r_{2})$-- atoms $\{ a_{Q} \}$ 
and a sequence of coefficients $\{ c'(Q) \} \in 
a^s(e^{s'}_{pq})^{\sigma}_{x_{0}}$ such that 
$f=\sum_{l(Q) \leq 1}c(Q)\phi_{Q}=\sum_{l(Q) \leq 1}c'(Q)a_{Q}$.
Hence, we see that $D \subset A$. 
  Therefore, we have 
 $A^s(\dot{E}^{s'}_{pq})^{\sigma}_{x_{0}} = M=A$.
We can prove (ii) by the same way in (i).
\qed

\bigskip

We recall the definition  of smooth wavelets. 

\bigskip

\noindent 
{\bf Definition 7.}\ \ Let $r \in \mathbb{N}_{0}$ and $L> n$.
A family of  $\{ \psi_{0},\ \psi^{(i)} \}$ is called ($r,\ L$)-- smooth wavelets if 
 $\{ \psi_{0}(x-k)\ (k \in \mathbb{Z}^n), \ 
2^{nj/2}\psi^{(i)}(2^{j}x-k)\ 
(i=1, \cdots, 2^n-1,  j \in {\mathbb N}_{0}, \ \
 k \in {\Bbb Z}^n) \}$ 
forms an orthonormal basis of 
$L^2(\mathbb{R}^n)$, and $\psi^{(i)}$ satisfies (5.1), (5.2) and (5.3), 
and  a scaling function  $\psi_{0}$ satisfies (5.4)

\bigskip

(5.1)\ \ $|\psi^{(i)}(x)| \leq C(1+|x|)^{-\max(L, L_{0})}$ for some $ L_{0}> n+r$,

(5.2)\ \ $|\partial^{\gamma}\psi^{(i)}(x)| \leq C(1+|x|)^{-L}$ 
for $0 <|\gamma| \leq r$,

(5.3)\ \ $\int_{\mathbb{R}^n} \psi^{(i)}(x) x^{\gamma} dx =0$ for $|\gamma| <r$

\noindent
where (5.2) and (5.3) are void when $r=0$.

(5.4)\ \ $|\partial^{\gamma}\psi_{0}(x)| \leq C(1+|x|)^{-L}$ 
for\ \  $|\gamma| \leq r$,

\noindent
but  $\psi_{0}$ does not satisfy the vanishing moment condition (5.3).
We will forget to write the index $i$ of the wavelet, 
which is of no consequence.

We put $\psi_{0,k}(x)= \psi_{0}(x-k),\ k \in \mathbb{Z}^n$,
$\psi_{Q}(x)= \psi(l(Q)^{-1}(x-x_{Q}))$ for a dyadic cube $Q$ with $l(Q) \leq 1$.

\bigskip

\noindent
{\bf Theorem 3.}\ \  {\it 
Let $s, \ s',\ \sigma \in {\Bbb R}, \ x_{0} \in \mathbb{R}^n$ and 
$0 < p,\ q \leq \infty$. 

{\rm (i)}\ \ For a family of ($r, L$)-- smooth wavelets $\{ \psi_{0},\ \psi \}$ satisfying 

\bigskip

{\rm (6.1)}\ \  $r > \max(s',\ \sigma+s + s'- \frac{n}{p},\  J-n-s')$
 and

{\rm (6.2)}\ \ $L> J$,  
\ \ 
where $J$ as in Lemma {\rm 1}, 

we have

\bigskip

$A^s(E^{s'}_{pq})^{\sigma}_{x_{0}}= \{ 
f = \sum_{k \in \mathbb{Z}^n}c_{k}\psi_{0,k}+\sum_{l(Q) \leq 1}c(Q)\psi_{Q} 
:  (c_{k}) \in a^s(e^{s'}_{pq})^{\sigma}_{x_{0}
}$,

$(c(Q)) \in a^s(e^{s'}_{pq})^{\sigma}_{x_{0}} \}$,

\bigskip

\noindent
where  $(c_{k})_{k \in \mathbb{Z}^n} \in a^s(e^{s'}_{pq})^{\sigma}_{x_{0}}$ means that 
$(c^{0}(Q))_{l(Q)\leq 1} \in a^s(e^{s'}_{pq})^{\sigma}_{x_{0}}$
such as $c^{0}(Q) =c_{k}$ if $Q=Q_{0,k}=[0, 1)^n+k,\ k \in \mathbb{Z}^n$ and 
$c^{0}(Q)=0$ if $l(Q) < 1$.

{\rm (ii)}\ \ For a family of ($r, L$)-- smooth wavelets  $\{ \psi_{0},\ \psi \}$ satisfying

\bigskip

{\rm $(6.1)'$}\ \  $r > \max(s'+(\sigma\vee 0),\ (\sigma\vee 0)+s + s'- \frac{n}{p},\  J-n-s'-(\sigma\wedge 0))$
 and

{\rm $(6.2)'$}\ \ $L> J+\sigma $

we have

\bigskip

$A^s(\tilde{E}^{s'}_{pq})^{\sigma}_{x_{0}}= \{ 
f = \sum_{k \in \mathbb{Z}^n}c_{k}\psi_{0,k}+\sum_{l(Q) \leq 1}c(Q)\psi_{Q} 
:  (c_{k}) \in a^s(\tilde{e}^{s'}_{pq})^{\sigma}_{x_{0}}$,

$(c(Q)) \in a^s(\tilde{e}^{s'}_{pq})^{\sigma}_{x_{0}} \}$.
}

\bigskip

\noindent
{\bf Remark 6.}\ \   We see that  by Lemma 3, 
$\sum_{k \in \mathbb{Z}^n}c_{k}\psi_{0,k}$  and 

\noindent
$\sum_{l(Q) \leq 1}c(Q)\psi_{Q}$ are convergent in $\mathcal{S}'$ for 
$(c_{k}),\ (c(Q)) \in a^s(e^{s'}_{pq})^{\sigma}_{x_{0}}$ or $a^s(\tilde{e}^{s'}_{pq})^{\sigma}_{x_{0}}$.

\bigskip

\noindent
{\it Proof of Theorem {\rm 3}.}\ \  (i)\ \ We may assume $\sigma \geq 0$ by Remark 1. We put 
$W= \{ 
f = \sum_{k \in \mathbb{Z}^n}c_{k}\psi_{0,k}+\sum_{l(Q) \leq 1}c(Q)\psi_{Q} 
:  (c_{k}), (c(Q)) \in a^s(e^{s'}_{pq})^{\sigma}_{x_{0}}
\}$.

Let $W \ni f = \sum_{k \in \mathbb{Z}^n}c_{k}\psi_{0,k}+\sum_{l(Q) \leq 1}c(Q)\psi_{Q}$ 
and  we consider the $\varphi$--transform

$$
\psi_{0,k}=\sum_{l(P)\leq 1} l(P)^{-n}\langle \psi_{0,k}\ ,\ \varphi_{P} \rangle \phi_{P}
$$
$$
\psi_{Q}=\sum_{l(P)\leq 1} l(P)^{-n}\langle \psi_{Q}\ ,\ \varphi_{P} \rangle \phi_{P}
$$
\noindent
where $\phi_{P}$ and $\varphi_{P}$ as in Remark 2.
Then we have

$$
f=\sum_{l(P) \leq 1} (B_{1}c_{k})(P) \phi_{P} + \sum_{l(P)\leq 1} (A_{1}c)(P) 
\phi_{P}
$$
where $B_{1}=\{ l(P)^{-n}\langle \psi_{0,k}\ ,\ \varphi_{P} \rangle \}_{Pk}$ and
\noindent
$A_{1}=\{ l(P)^{-n}\langle \psi_{Q}\ ,\ \varphi_{P} \rangle \}_{PQ}$.
Lemma 1 and Lemma 2 yield that 
$B_{1}$ and $A_{1}$ are almost diagonal and 
$B_{1}c_{k},\ A_{1}c \in a^s(e^{s'}_{pq})^{\sigma}_{x_{0}}$
for $c_{k},\ c \in a^s(e^{s'}_{pq})^{\sigma}_{x_{0}}$.
Hence, 
by  Theorem 1, we see that 
$W \subset D =
 A^s(E^{s'}_{pq})^{\sigma}_{x_{0}}$ where $D$ is as 
in the proof of Theorem 1.

Conversely, 
let $D \ni f = \sum_{l(Q)\leq 1} c(Q) \phi_{Q}$ and we consider the 
wavelet expansion

$$
\phi_{Q}=\sum_{k \in \mathbb{Z}^n} \langle \phi_{Q}\ ,\ \psi_{0,k} \rangle \psi_{0,k}+
\sum_{l(P)\leq 1} l(P)^{-n}\langle \phi_{Q}\ ,\ \psi_{P} \rangle \psi_{P}.$$

\noindent
Then we have

$$
f=\sum_{k \in \mathbb{Z}^n} (B_{2}c)(k) \psi_{0,k}+
\sum_{l(Q)\leq 1} (A_{2}c)(Q) \phi_{Q} $$
where $B_{2}=\{ \langle \phi_{Q}\ ,\ \psi_{0,k} \rangle \}_{kQ}$ and 
$A_{2}=\{ l(P)^{-n}\langle \phi_{Q}\ ,\ \psi_{P} \rangle \}_{PQ}$.
Lemma 1 and Lemma 2 yield that $B_{2}$ and $A_{2}$ are 
almost diagonal and 
$B_{2}c\ A_{2}c \in a^s(e^{s'}_{pq})^{\sigma}_{x_{0}}$
for $c \in a^s(e^{s'}_{pq})^{\sigma}_{x_{0}}$.
Hence, by  Theorem 1, we see that 
$A^s(\dot{E}^{s'}_{pq})^{\sigma}_{x_{0}}=D \subset W$.  

We can prove (ii) by the same way in (i). 
Hence we obtain  the result of Theorem 3.
\qed

\bigskip

\noindent
{\bf Remark 7.}\ \  
(1)\ \ we see that Theorem 3 is independent of the choice of 
  smooth wavelets $\{ \psi_{0},\ \psi^{(i)} \}$ (see Remark 3 (2)).   
 
(2)\ \
For $f \in A^s(E^{s'}_{pq})^{\sigma}_{x_{0}}$ or $A^s(\tilde{E}^{s'}_{pq})^{\sigma}_{x_{0}}$ the parings $
\langle f\ ,\ \psi_{0, k}\rangle$ and $\langle f\ ,\ \psi_{Q} \rangle$ are well-defined. More explicitly, we see that for any $\{\phi_{Q},\ \ \varphi_{Q} \}$ as in Remark 2, 
$$
\langle f\ ,\ \psi_{0, k}\rangle= \sum_{l(P) \leq 1}l(P)^{-n}
\langle f\ ,\ \phi_{P}\rangle \langle \psi_{0, k}\ ,\ \varphi_{P} \rangle
\equiv \sum_{l(P) \leq 1}c(f)(P)\langle \psi_{0,k}\ ,\ \varphi_{P} \rangle
$$
and
$$ 
\langle f\ ,\ \psi_{Q} \rangle=\sum_{l(P) \leq 1}l(P)^{-n}
\langle f\ ,\ \phi_{P} \rangle\langle \psi_{Q}\ ,\ \varphi_{P} \rangle\equiv \sum_{l(P) \leq 1}c(f)(P)
 \langle \psi_{Q}\ ,\ \varphi_{P} \rangle
$$
are convergent by  Lemma 3 and (b) in the proof of Theorem 1. 
Thus, for $f \in A^s(E^{s'}_{pq})^{\sigma}_{x_{0}}$ or $A^s(\tilde{E}^{s'}_{pq})^{\sigma}_{x_{0}}$ 
 we have a wavelet expansion 
$f= \sum_{k \in \mathbb{Z}^n}c_{k}\psi_{0,k}+
\sum_{l(Q)\leq 1}c(Q)\psi_{Q}$ in $\mathcal{S}'$ 
and its representation  is unique in $\mathcal{S}'$, that is,  
$c_{k}=\langle f\ ,\ \psi_{0, k}\rangle$ and 
$c(Q)=l(Q)^{-n}\langle f\ ,\ \psi_{Q} \rangle$.
 Hence, we have that  by Lemma 1, Lemma 2 and (b) in the proof of Theorem 1,
\begin{eqnarray*}
\lefteqn{
||(c_{k})||_{a^s(e^{s'}_{pq})^{\sigma}_{x_{0}}}=
||\langle f\ ,\ \psi_{0, k}\rangle||_{a^s(e^{s'}_{pq})^{\sigma}_{x_{0}}}
}
\\
&\leq&
 ||\sum_{l(P) \leq 1}c(f)(P)\langle \psi_{0,k}\ ,\ \varphi_{P} \rangle||_{a^s(e^{s'}_{pq})^{\sigma}_{x_{0}}}
\\
&\leq&
 C||c(f)||_{a^s(e^{s'}_{pq})^{\sigma}_{x_{0}}}
\leq C||f||_{A^s(E^{s'}_{pq})^{\sigma}_{x_{0}}}
\end{eqnarray*}
and
\begin{eqnarray*}
\lefteqn{
||(c(Q))||_{a^s(e^{s'}_{pq})^{\sigma}_{x_{0}}}=||l(Q)^{-n}
\langle f\ ,\ \psi_{Q} \rangle||_{a^s(e^{s'}_{pq})^{\sigma}_{x_{0}}}
}
\\
&\leq&
 C||\sum_{l(P) \leq 1}c(f)(P)l(Q)^{-n}\langle \psi_{Q}\ ,\ \varphi_{P} \rangle||_{a^s(e^{s'}_{pq})^{\sigma}_{x_{0}}}
\\
&\leq&
 C||c(f)||_{a^s(e^{s'}_{pq})^{\sigma}_{x_{0}}}
\leq C||f||_{A^s(E^{s'}_{pq})^{\sigma}_{x_{0}}}.
\end{eqnarray*}
Conversely, we consider the $\varphi$-transform

$$
\psi_{0,k}=\sum_{P} l(P)^{-n}\langle \psi_{0,k}\ ,\ \varphi_{P} \rangle \phi_{P}$$
and
$$
\psi_{Q}=\sum_{P} l(P)^{-n}\langle \psi_{Q}\ ,\ \varphi_{P} \rangle \phi_{P}.
$$ 
Then we have 
$$
f= \sum_{k \in \mathbb{Z}^n}c_{k}\psi_{0,k}+
\sum_{Q}c(Q)\psi_{Q}
= \sum_{k \in \mathbb{Z}^n}(Bc_{k})(P)\phi_{P}+
\sum_{Q}Ac(P)\phi_{P}
$$ 
where $B=\{l(P)^{-n}\langle \psi_{0,k}\ ,\ \varphi_{P} \rangle \}$ and 
$A=\{l(P)^{-n}\langle \psi_{Q}\ ,\ \varphi_{P} \rangle \}$.
Hence we have by Lemma 1, Lemma 2 and (a) in the proof of Theorem 1,
\begin{eqnarray*}
\lefteqn{
||f||_{A^s(E^{s'}_{pq})^{\sigma}_{x_{0}}}\leq C||(Bc_{k})+(Ac)||_{a^s(e^{s'}_{pq})^{\sigma}_{x_{0}}}
}
\\
&\leq&
C||Bc_{k}||_{a^s(e^{s'}_{pq})^{\sigma}_{x_{0}}}+C||Ac||_{a^s(e^{s'}_{pq})^{\sigma}_{x_{0}}}
\\
&\leq&
 C||c_{k}||_{a^s(e^{s'}_{pq})^{\sigma}_{x_{0}}}+C||c||_{a^s(e^{s'}_{pq})^{\sigma}_{x_{0}}}.
\end{eqnarray*}
Therefore, we have
$$||f||_{A^s(E^{s'}_{pq})^{\sigma}_{x_{0}}} \sim
||(c_{k})||_{a^s(e^{s'}_{pq})^{\sigma}_{x_{0}}}
+||(c(Q))||_{a^s(e^{s'}_{pq})^{\sigma}_{x_{0}}}.
$$ 
Similarly, we also obtain 
$$||f||_{A^s(\tilde{E}^{s'}_{pq})^{\sigma}_{x_{0}}} \sim
||(c_{k})||_{a^s(\tilde{e}^{s'}_{pq})^{\sigma}_{x_{0}}}
+||(c(Q))||_{a^s(\tilde{e}^{s'}_{pq})^{\sigma}_{x_{0}}}.
$$ 

\section{Applications}

\bigskip

\noindent
{\bf Definition 8.}  
Let $\mathcal{T}$ be the space of Schwartz test functions ($C^{\infty}$-functions with compact support) 
and $\mathcal{T}'$ its dual.  
For arbitrary $r_{1},\ r_{2} \in {\mathbb N}_{0}$ 
the Calder$\acute{{\rm o}}$n--Zygmund operator $T$ with an exponent $\epsilon 
> 0$  is 
a continuous linear operator $\mathcal{T} \rightarrow \mathcal{T}'$ such that 
its kernel $K$ off the diagonal $\{ (x,y) \in {\mathbb R}^n \times 
{\mathbb R}^n : x=y \}$ satisfies

(7.1)\ \ $|\partial^{\gamma}_{1}K(x,y)| \leq C|x-y|^{-(n+|\gamma|)}$ 
for  $|\gamma| \leq  r_{1}$,

(7.2)\ \ $|K(x,\ y)-K(x,\ y')|
\leq  C|y-y'|^{r_{2} +\epsilon}|x-y|^{-(n+r_{2} + \epsilon)}$ 
if $2|y'-y| \leq |x-y|$,

(7.3)\ \ $|\partial^{\gamma}_{1}K(x,\ y)-\partial^{\gamma}_{1}K(x,\ y')|
\leq  C|y-y'|^{\epsilon}|x-y|^{-(n+|\gamma| + \epsilon)}$ 

if $2|y'-y| \leq |x-y|$ for $0 <|\gamma| \leq r_{1}$ 

(where this statement is void when $r_{1}=0$),

 $|\partial^{\gamma}_{1}K(x,\ y)-\partial^{\gamma}_{1}K(x',\ y)| 
\leq  C|x'-x|^{\epsilon}|x-y|^{-(n+|\gamma| + \epsilon)}$ 

if $2|x'-x| \leq |x-y|$ for $|\gamma| \leq r_{1}$,

\noindent
(where the subindex 1 stands for derivatives in the first variable)

(7.4)\ \ $T$ is bounded on $L^2({\Bbb R}^n)$.

\bigskip

We obtain the following theorem.

\bigskip
\noindent
{\bf Theorem 4.} \ \ {\it 
 Let $s,\ s',\ \sigma  \in {\Bbb R}, \   x_{0} \in \mathbb{R}^n$,
   $0< p,\ q \leq \infty$, $r_{1},\ r_{2} \in {\mathbb N}_{0}$ and 
 $J$ as in Lemma {\rm 1}.  

{\rm (i)}\ \ 
 The Calder$\acute{o}$n--Zygmund operator $T$ with an exponent $\epsilon 
>J -n$  
 satisfying $T(x^{\gamma})=0$ for $|\gamma| \leq r_{1}$ and 
  $T^{*}(x^{\gamma})=0$ for $|\gamma| < r_{2}$,  
 is bounded on  
$A^s(E^{s'}_{pq})^{\sigma}_{x_{0}}$
if $r_{1}$ and $r_{2}$ satisfy {\rm (4.1)} and {\rm (4.2)} as in Theorem {\rm 2} respectively.
 
{\rm (ii)}\  \ The Calder$\acute{o}$n--Zygmund operator $T$ with an exponent $\epsilon 
>J -n+\sigma$  
 satisfying $T(x^{\gamma})=0$ for $|\gamma| \leq r_{1}$ and 
  $T^{*}(x^{\gamma})=0$ for $|\gamma| < r_{2}$,  
 is bounded on  
$A^s(\tilde{E}^{s'}_{pq})^{\sigma}_{x_{0}}$ 
if $r_{1}$ and $r_{2}$ satisfy {\rm (4.1)'} and {\rm (4.2)'} as in Theorem {\rm 2} respectively.
}

\bigskip

\noindent
{\it Proof}.\ \  
The proof is similar to ones of [12]. 
 
(i)\ \ We may assume $\sigma \geq 0$ by Remark 1. 
Let $f \in A^s(E^{s'}_{pq})^{\sigma}_{x_{0}}$. Then 
 we consider a wavelet expansion 
$f = \sum_{k}c_{k}\psi_{0,k}+ \sum_{l(Q)\leq 1}c(Q)\psi_{Q}: \  (c_{k}),\ (c(Q)) \in a^s(e^{s'}_{pq})^{\sigma}_{x_{0}}$ from Theorem 3. 
We may suppose that smooth wavelets $\{ \psi_{0}\  \psi\}$ are 
compactly supported by Remark 7 (1). 
Then there exists a positive constant c such that  
supp $\psi_{0,k} \subset c Q_{0,k}$ where $Q_{0,k}=[0, 1)^n+k$ and supp $\psi_{Q} \subset cQ$ for every dyadic cube $Q$ with $l(Q)=2^{-l} \leq 1$.  

We claim that $Tf = \sum_{k}c_{k}(T\psi_{0,k})+ \sum_{l(Q)\leq 1}c(Q)(T\psi_{Q}) \equiv \sum_{k}c_{k}m_{k}+\sum_{l(Q) \leq 1}c(Q)m_{Q}$ is convergent in 
$\mathcal{S}'$ and  
$||Tf||_{A^s(E^{s'}_{pq})^{\sigma}_{x_{0}}} \leq  
C||f||_{A^s(E^{s'}_{pq})^{\sigma}_{x_{0}}}$. 

More precisely , we will show that $m_{k}$ and 
$m_{Q}$ satisfy following properties: 

(8.1)\ \ 
$
|m_{k}(x)| \leq C(1+l(Q)^{-1}|x-x_{k}|)^{-L}
$
with $L>J$,

(8.2)\ \ 
$
|m_{Q}(x)| \leq C(1+l(Q)^{-1}|x-x_{Q}|)^{-(n+r_{2}+\epsilon)}
$,

(8.3)\ \ 
$|\partial^{\gamma}m_{Q}(x)| \leq Cl(Q)^{-|\gamma|}(1+l(Q)^{-1}|x-x_{Q}|)
^{-(n+\epsilon)}
$
 for\ \  $0 <|\gamma| \leq r_{1}$, and

(8.4)\ \ $\int_{{\mathbb R}^n}  x^{\gamma}m_{Q}(x) dx =0$ for\ \  $|\gamma| < r_{2}$.

From the assumption $T^{*}x^{\gamma}=0$ 
for $|\gamma| < r_{2}$ we have 
$\int_{{\mathbb R}^n}  x^{\gamma}m_{Q}(x) dx =0$ for\ \  $|\gamma| < r_{2}$, 
that is, (8.4) holds.

 We choose a suitable large constant $C_{0}$. 
From Fraizer--Torres--Weiss [12: Corollary 2.14],  when $|x-x_{Q}| < 2C_{0}2^{-l}$, we have 
$$|m_{k}(x)|\leq 
||m_{k}||_{\infty} \leq 
C\sum_{|\beta|\leq 1}||\partial^{\beta}\psi_{0,k}||_{\infty} 
\leq  C  \leq C (1+|x-x_{k}|)^{-L}
$$ 
and
$$
|\partial^{\gamma}m_{Q}(x)|\leq 
||\partial^{\gamma}m_{Q}||_{\infty} \leq 
C\sum_{|\alpha|\leq |\gamma|+1}2^{l(|\gamma|-|\alpha|)}2^{l|\alpha|}
||\partial^{\alpha}\psi_{Q}||_{\infty} $$
$$
\leq 
C 2^{l|\gamma|} \leq 
C l(Q)^{-|\gamma|}(1+l(Q)^{-1}|x-x_{Q}|)^{-L}
$$ 
for any $L \geq 0$ and $|\gamma| \leq r_{1}$.
 When $|x-x_{Q}| \geq 2C_{0}2^{-l}$, using  (7.1) and (7.2) in Definition 8, 
we obtain 
\begin{eqnarray*}
\lefteqn{
|m_{k}(x)| = |\int_{\mathbb{R}^n}K(x,y)\psi_{0,k}(y) dy|
\leq C\int_{\mathbb{R}^n}|K(x,y)||\psi_{0,k}(y)| dy 
} 
\\
&\leq&
C\int_{|y-x_{k}| \leq C_{0}}|x-y|^{-n}(1+|y-x_{k}|)^{-L} dy
\leq
C(1+|x-x_{k}|)^{-(L+n)}.
\end{eqnarray*}
Moreover, using  (7.3) in Definition 8
for $0 < |\gamma| \leq r_{1}$, 
we have 
\begin{eqnarray*}
\lefteqn{
|\partial^{\gamma}m_{Q}(x)| \leq 
C\int_{|y-x_{Q}| \leq C_{0}2^{-l}}
|\partial^{\gamma}_{1}K(x,y)-\partial^{\gamma}_{1}K(x,x_{Q})||\psi_{Q}(y)| dy  
}
\\
&\leq&
 C\int_{|y-x_{Q}| \leq C_{0}2^{-l}}|y-x_{Q}|^{\epsilon}
|x-x_{Q}|^{-(n+|\gamma|+\epsilon)} dy 
\\
&\leq& 
C2^{-l(n+\epsilon)}|x-x_{Q}|^{-(n+|\gamma|+\epsilon)} 
\leq
C2^{l|\gamma|}(1+2^{l}|x-x_{Q}|)^{-(n+\epsilon)}.
\end{eqnarray*}
Therefore, we obtain (8.1), (8,2) and (8.3).
Hence by  Lemma 3, $Tf =\sum_{k}c_{k}m_{k}+\sum_{Q}c(Q)m_{Q}$ is convergent in 
$\mathcal{S}'$ from (8.1), (8.2), (8.3) and (8.4). 
For the wavelet expansion

$$
m_{k}=\sum_{k}\langle m_{k}\ ,\ \psi_{0,k} \rangle \psi_{0,k}+\sum_{P} l(P)^{-n}\langle m_{k}\ ,\ \psi_{P} \rangle \psi_{P},
$$

$$
m_{Q}=\sum_{k} \langle m_{Q}\ ,\ \psi_{0,k} \rangle \psi_{0,k}+\sum_{P} l(P)^{-n}\langle m_{Q}\ ,\ \psi_{P} \rangle \psi_{P},
$$
we have 
\begin{eqnarray*}
\lefteqn{
Tf=\sum_{k} c_{k} m_{k}+\sum_{l(Q)\leq 1} c(Q) m_{Q} = 
}
\\
&&
\sum_{k} ((B_{1}c_{k})+(B_{2}c_{k})) \psi_{0,k}+\sum_{l(P)\leq 1} ((A_{1}c)+(A_{2}c))(P) \psi_{P}
\end{eqnarray*}
where $B_{1}=\{ \langle m_{k}\ ,\ \psi_{0,k'} \rangle \}_{k'k},\ 
B_{2}=\{ \langle m_{Q}\ ,\ \psi_{0,k'} \rangle \}_{k'Q},
$

\noindent
$
\ A_{1}=
\{ l(P)^{-n}\langle m_{k}\ ,\ \psi_{P} \rangle \}_{Pk},\ A_{2}=
\{ l(P)^{-n}\langle m_{Q}\ ,\ \psi_{P} \rangle \}_{PQ}$.
By Lemma 1 , Lemma 2, (8.1), (8,2), (8,3) and (8.4)  the operators $B_{1},\ B_{2},\ 
A_{1},\ A_{2}$ are bounded on $a^s(e^{s'}_{pq})^{\sigma}_{x_{0}}$ if $r_{1}$ and $r_{2}$ satisfy (4.1) and (4.2) respectively. 
By Remark 7 (2), 
it follows that 

$$||Tf||_{A^s(E^{s'}_{pq})^{\sigma}_{x_{0}}} \sim 
||(B_{1}c_{k}+B_{2}c_{k})||_{a^s(e^{s'}_{pq})^{\sigma}_{x_{0}}}
+||(A_{1}c+A_{2}c)||_{a^s(e^{s'}_{pq})^{\sigma}_{x_{0}}}
$$
$$
\leq C(||c_{k}||_{a^s(e^{s'}_{pq})^{\sigma}_{x_{0}}}+ ||c||_{a^s(e^{s'}_{pq})^{\sigma}_{x_{0}}}) \sim 
C||f||_{A^s(E^{s'}_{pq})^{\sigma}_{x_{0}}}.
$$ 
Similarly, we obtain (ii).  These complete the proof.
\qed

\bigskip

\noindent
{\bf Definition 9.}\ \ 
Let $\mu \in \mathbb{R}$. A smooth function $a$ defined on $\mathbb{R}^n \times \mathbb{R}^n$ is said to belong to the class $S^{\mu}_{1,1}(\mathbb{R}^n)$ if $a$ satisfies the following differential inequalities that for all $\alpha,\ \beta \in \mathbb{N}_{0}^n$,
$$
\sup_{x, \xi}(1+|\xi|)^{-\mu-|\alpha|+|\beta|}|\partial^{\alpha}_{x}
\partial^{\beta}_{\xi}a(x,\ \xi)|< \infty.
$$
$a(x, D)$ is the corresponding pseudo-differential operator such that 
$$
a(x, D)f(x)=\int_{\mathbb{R}^n}e^{ix\xi}a(x, \xi)\hat{f}(\xi)\ d\xi
$$
for  $f \in \mathcal{S}$.

\bigskip

\noindent
{\bf Theorem 5}. \ \ {\it 
 Let $s,\ s',\ \sigma \in {\Bbb R},\   x_{0} \in \mathbb{R}^n$,
   $0< p,\ q \leq \infty$.
Let $\mu \in \mathbb{R}$, $J$  as in Lemma {\rm 1} and $a \in S^{\mu}_{1,1}(\mathbb{R}^n)$.

{\rm (i)}\ \ 
$a(x,D)$ is a continuous linear mapping 
from 
$A^s(E^{s'}_{pq})^{\sigma}_{x_{0}}$ to $A^s(E^{s'-\mu}_{pq})^{\sigma}_{x_{0}}$
 if $s' > J-n+\mu$  or  $a(x,\xi)=a(\xi)$. 

{\rm (ii)}\ \ 
$a(x,D)$ is a continuous linear mapping 
from 
$A^s(\tilde{E}^{s'}_{pq})^{\sigma}_{x_{0}}$ to $A^s(\tilde{E}^{s'-\mu}_{pq})^{\sigma}_{x_{0}}$ 
 if $s' > J-n+\mu+\sigma \wedge 0$  or  $a(x,\xi)=a(\xi)$. 
}

\bigskip

\noindent
{\it Proof}.\ \  (i)\ \ We may assume $\sigma \geq 0$ by Remark 1. 
We write $T\equiv a(x, D)$. Let $f \in 
A^s(E^{s'}_{pq})^{\sigma}_{x_{0}}$. By Theorem 1, we consider the $\varphi$-transform 
$f = \sum_{P}c(P)\phi_{P}$  where $c(P)=c(f)(P)=l(P)^{-n}\langle f,\  \varphi_{P} \rangle$ and $\phi_{P},\ \varphi_{P}$ as in Remark 2. Then we see $(c(P)) 
\in a^s(e^{s'}_{pq})^{\sigma}_{x_{0}}$. 
We write  that $Tf= \sum_{P}c(P)m_{P}$ where $m_{P}=T\phi_{P}$. 
We see  for a dyadic cube $P$ with $l(P)=2^{-j}$

$$
m_{P}=\int e^{ix\xi}a(x,\xi)\hat{\phi}_{P}(\xi)\ d\xi.
$$
Then we have, using a change of variables,  
$$
m_{P}(x)=\int e^{i(x-x_{P})(2^j\xi)}a(x, 2^j\xi)\hat{\phi}(\xi)\ d\xi.
$$ 
By the fact that  $(1-\triangle_{\xi})^L(e^{ix\xi})=(1+|x|^2)^Le^{ix\xi}$ for 
the Laplacian $\triangle$ and using an integration by parts, we obtain for $\gamma \in \mathbb{N}_{0}^n$ and $l(P) <1$,
\begin{eqnarray*}
\lefteqn{\partial_{x}^{\gamma}m_{P}(x)
} 
\\
&=&
\int (1-\triangle_{\xi})^L(e^{i2^j(x-x_{P})\xi})(1+(2^j|x-x_{P}|)^2)^{-L}\times
\\
&&
\sum_{\delta \leq \gamma}(2^ji\xi)^{\delta}\partial_{x}^{\gamma-\delta}a(x, 2^j\xi)\hat{\phi}(\xi)\ d\xi 
\\
&=&
C(1+(2^j|x-x_{P}|)^2)^{-L}\int e^{i2^j(x-x_{P})\xi}
(1-\triangle_{\xi})^L \times
\\
&&
\sum_{\delta \leq \gamma}(2^ji\xi)^{\delta}\partial_{x}^{\gamma-\delta}a(x, 2^j\xi)
\hat{\phi}(\xi)\ d\xi. 
\end{eqnarray*}
Thus, we have 
\begin{eqnarray*}
\lefteqn{|\partial_{x}^{\gamma}m_{P}(x)|
} 
\\
&\leq&
C(1+2^j|x-x_{P}|)^{-2L}\int 
\sum_{\delta \leq \gamma}\sum_{|\alpha+\beta+\tau|\leq 2L, \alpha \leq \delta }
\times
\\
&&
2^{j|\delta|}2^{j|\beta|}
|\partial_{\xi}^{\alpha}(\xi)^{\delta}||\partial_{\xi}^{\beta}\partial_{x}^{\gamma-\delta}a(x, 2^j\xi)||\partial_{\xi}^{\tau}\hat{\phi}(\xi)|\ d\xi 
\\
&\leq&
C(1+2^j|x-x_{P}|)^{-2L}
\int 
\sum_{\delta \leq \gamma}\sum_{|\alpha+\beta+\tau|\leq 2L, \alpha \leq \delta}
\times
\\
&&
2^{j|\delta|}2^{j|\beta|}
|\xi|^{|\delta|-|\alpha|}(1+2^j|\xi|)^{\mu+|\gamma|-|\delta|-|\beta|}|\partial_{\xi}^{\tau}\hat{\phi}(\xi)|\ d\xi 
\\
&\leq&
C2^{j\mu}2^{j|\gamma|}(1+2^j|x-x_{P}|)^{-2L}
\end{eqnarray*}
and  similarly, for $P$ with $l(P)=1$,
\begin{eqnarray*}
\lefteqn{|\partial_{x}^{\gamma}m_{P}(x)|
} 
\\
&\leq&
C(1+|x-x_{P}|)^{-2L}\times
\\
&&
\int 
\sum_{|\alpha+\beta+\tau|\leq 2L}
(1+|\xi|)^{\mu+|\gamma|-|\alpha|-|\beta|}|\partial_{\xi}^{\tau}\hat{\phi_{0}}(\xi)|\ d\xi 
\\
&\leq&
C(1+|x-x_{P}|)^{-2L}.
\end{eqnarray*}

Hence, $m_{P}(x)$ satisfies
$$
|2^{-j\mu}\partial^{\gamma}m_{P}(x)| \leq C2^{j|\gamma|}(1+2^j|x-x_{P}|)^{-2L}$$
for $P$ with $l(P) \leq 1,\ $any $\gamma \in \mathbb{N}_{0}$ and any $L \geq 0$. 
We choose a suitable large $L$. 
For the $\varphi$-- transform
$$
2^{-j\mu}m_{P}=\sum_{l(R) \leq 1} l(R)^{-n}\langle 2^{-j\mu}m_{P}\ ,\ \varphi_{R} \rangle \phi_{R},
$$
we have
$$
Tf=\sum_{l(P) \leq 1}2^{j\mu} c(P) (2^{-j\mu}m_{P}) = \sum_{l(R) \leq 1} A(2^{j\mu}c)(R) \phi_{R},
$$
where $A=\{ l(R)^{-n}\langle 2^{-j\mu}m_{P}\ ,\ \varphi_{R} \rangle \}_{RP}$.
From Lemma 1 and Lemma 2, 
$A$ is bounded on $a^s(e^{s'-\mu}_{pq})^{\sigma}_{x_{0}}$ if $s' > J-n+\mu$ or  $a(x,\xi)=a(\xi)$.  We remark that in the case $s' > J-n+\mu$, we do not assume the vanishing moment condition for $m_{P}$. 
But  in the case  $a(x,\xi)=a(\xi)$, we have the vanishing moment condition for $m_{P}$, indeed,   for any $P$ with $l(P) < 1$,
$\int x^{\gamma}m_{P}(x)\ dx=C\partial^{\gamma}\hat{m}_{P}(0)=C\partial^{\gamma}(\hat{\phi}_{P}) \cdot a)(0)=0$ for any $\gamma \in \mathbb{N}_{0}$. 
From (a) and (b) in the proof of Theorem 1, it follows that 
\begin{eqnarray*}
\lefteqn{
||Tf||_{A^s(E^{s'-\mu}_{pq})^{\sigma}_{x_{0}}} \leq 
C||A(2^{j\mu}c)||_{a^s(e^{s'-\mu}_{pq})^{\sigma}_{x_{0}}}
}
\\
&\leq&
 C||2^{j\mu}c||_{a^s(e^{s'-\mu}_{pq})^{\sigma}_{x_{0}}}
\leq C||c||_{a^s(e^{s'}_{pq})^{\sigma}_{x_{0}}}
\leq
C||f||_{A^s(E^{s'}_{pq})^{\sigma}_{x_{0}}}.
\end{eqnarray*}

(ii)\ \ Similarly, we can prove for this case. 
\qed

\bigskip

\noindent
{\bf Corollary }.\ \ {\it
Let $s,\ s',\ \sigma \in {\Bbb R},  \ x_{0} \in \mathbb{R}^n$,  $0< p,\ q \leq \infty$.

{\rm (i)}\ \  Let $\mu \in \mathbb{R}$. Then the Bessel potential 
 $(1-\triangle)^{\mu/2}$ is a continuous isomorphisms 
 from $A^s(E^{s'}_{pq})^{\sigma}_{x_{0}}$ onto $A^s(E^{s'-\mu}_{pq})^{\sigma}_{x_{0}}$, 
 and from 
$A^s(\tilde{E}^{s'}_{pq})^{\sigma}_{x_{0}}$ onto $A^s(\tilde{E}^{s'-\mu}_{pq})^{\sigma}_{x_{0}}$. 

{\rm (ii)} \ \ 
Let $\gamma \in \mathbb{N}^n_{0}$. 
Then the differential operator $\partial^{\gamma}$ 
 is  continuous  
 from $A^s(E^{s'}_{pq})^{\sigma}_{x_{0}}$ to $A^s(E^{s'-|\gamma|}_{pq})^{\sigma}_{x_{0}}$, 
 and from 
$A^s(\tilde{E}^{s'}_{pq})^{\sigma}_{x_{0}}$ to $A^s(\tilde{E}^{s'-|\gamma|}_{pq})^{\sigma}_{x_{0}}$. 
}

\bigskip

\noindent
{\it Proof}.\ \  These are immediate corollaries of Theorem 5.  
 To finish the proof of (i) we need to show the mapping is 
surjective and one to one.   For 
$h \in  A^s(E^{s'-\mu}_{pq})^{\sigma}_{x_{0}}$, we set $f = (1-\triangle)^{-\mu/2}h$. Then $h = (1-\triangle)^{\mu/2}f$. 
 \qed

\section{Characterizations via differences and oscillations}
\noindent
{\bf Definition 10.}
Let $k \in \mathbb{N}_{0}$. 
We define the differences of functions 

\bigskip

$
\triangle^1_{u}f(x) = f(x+u)-f(x)
$ and 
$
\triangle^{k+1}=\triangle^{1}\triangle^{k}.
$

\bigskip

We set 

\bigskip

$d^{k}_{i}f(y)=\frac{1}{|B_{i}(y)|}\int_{k|u| \leq 2^{-i}}
|\triangle^{k}_{u}f(y)|\ du$

\bigskip

\noindent
where $B_{i}(x)$ is the ball with a center $x$ and a radius $2^{-i}$, and 
$|B_{i}(x)|$ means its volume. 
It is obvious that $|d^{k}_{i}f(y)|\leq C \sup_{k|u| \leq 2^{-i}}
|\triangle^{k}_{u}f(y)|$.

We define the oscillation of   locally $L^p$ integrable functions $f$ 
($0 < p \leq \infty$) by

$$
 {\rm osc}^k_pf(x, i)=\inf (\frac{1}{|B_{i}(x)|}\int_{B_{i}(x)}|f(y)-P(y)|^p\ dy)^{1/p}
$$
with the suitable modification for $p=\infty$, 
where the infimum is taken  over all polynomials $P(x) \in \mathcal{P}_{k}$, 
the space of all polynomials with $\deg \leq k$ on ${\Bbb R}^n$.
By $P_{B}f$ for a ball $B$ we denote the unique polynomial in $\mathcal{P}_{k}$ such that 
$ \int_{B}(f(x)-P_{B}f(x))x^{\alpha}\ dx =0$
for all $|\alpha| \leq k$. We see that $||P_{B}f||_{L^{\infty}(B)} \leq \frac{1}{|B|}\int_{B}
|f(x)|\ dx $ and $P_{B}f =f$  for  $f \in \mathcal{P}_{k}$. 
 We put 

$$
 \Omega^k_pf(x, i)=(\frac{1}{|B_{i}(x)|}\int_{B_i(x)}|f(y)-P_{B_i(x)}f(y)|^p\ dy)^{1/p}.
$$

Then we see ${\rm osc}^k_pf(x, i)\sim \Omega^k_pf(x, i)$ if $1 \leq p \leq \infty$ (cf. [19]).

\bigskip

\noindent
{\bf Lemma 5.}\ \ {\it 
{\rm (i)}\ \  Let $s \in \mathbb{R},\ \sigma \geq 0$ and let $k \in \mathbb{N},\ \ k> s' > 0, \ \ 1 \leq p \leq \infty,\ \  0 < q \leq \infty$ and  let $f$ be locally $L^p$ integrable. 

Then we have 

\begin{eqnarray*}\lefteqn{ 
\sup_{ \mathcal{D} \ni Q \ni x_{0} }l(Q)^{-\sigma}
\sup_{ \mathcal{D} \ni P \subset 3Q\ }l(P)^{-s}\times
}
\\
&&
(\sum_{i \geq (-\log_{2}l(P))\vee 0}(2^{is'}\sup_{k|u| \leq 2^{-i}}||\triangle^{k}_uf||_{L^p(P)})^q)^{1/q} 
\\
&\leq&
 C
\sup_{ \mathcal{D} \ni Q \ni x_{0} }l(Q)^{-\sigma}
\sup_{ \mathcal{D} \ni P \subset 3Q\ }l(P)^{-s}\times
\\
&&
(\sum_{i \geq (-\log_{2}l(P))\vee 0}(2^{is'}
||{\rm osc}^{k-1}_pf(x, i)||_{L^p(P)})^q)^{1/q}, 
\end{eqnarray*}
and

\begin{eqnarray*}\lefteqn{
\sup_{ \mathcal{D} \ni Q \ni x_{0} }l(Q)^{-\sigma}
\sup_{ \mathcal{D} \ni P \subset 3Q\ }l(P)^{-s}
\times
}
\\
&&
(\sum_{i \geq (-\log_{2}l(P))\vee 0}(2^{is'}
||{\rm osc}^{k-1}_p
f(x, i)||_{L^p(P)})^q)^{1/q} 
\\
&\leq&
 C
\sup_{ \mathcal{D} \ni Q \ni x_{0} }l(Q)^{-\sigma}
\sup_{ \mathcal{D} \ni P \subset 3Q\ }l(P)^{-s}\times
\\
&&
(||f||_{L^p(P)}+(\sum_{i \geq (-\log_{2}l(P))\vee 0}(2^{is'}
||d^{k}_if||_{L^p(P)})^q)^{1/q}). 
\end{eqnarray*} 

{\rm (ii)}\ \  Let $s \in \mathbb{R},\ \sigma \geq 0$ and let $k \in \mathbb{N},\ \ k> s' > 0, \ \ 1 \leq p < \infty,\ \  0 < q \leq \infty$ and  let $f$ be locally $L^p$ integrable . Then we have 

\begin{eqnarray*}\lefteqn{ 
\sup_{ \mathcal{D} \ni Q \ni x_{0} }l(Q)^{-\sigma}
\sup_{ \mathcal{D} \ni P \subset 3Q\ }l(P)^{-s}\times
}
\\
&&
||(\sum_{i \geq (-\log_{2}l(P))\vee 0}(2^{is'}\sup_{k|u| \leq 2^{-i}}|\triangle^{k}_uf|)^q)^{1/q}||_{L^p(P)} 
\\
&\leq&
 C
\sup_{ \mathcal{D} \ni Q \ni x_{0} }l(Q)^{-\sigma}
\sup_{ \mathcal{D} \ni P \subset 3Q\ }l(P)^{-s}\times
\\
&&
||(\sum_{i \geq (-\log_{2}l(P))\vee 0}(2^{is'}
{\rm osc}^{k-1}_pf(x, i))^q)^{1/q}||_{L^p(P)}, 
\end{eqnarray*}
and

\begin{eqnarray*}\lefteqn{
\sup_{ \mathcal{D} \ni Q \ni x_{0} }l(Q)^{-\sigma}
\sup_{ \mathcal{D} \ni P \subset 3Q\ }l(P)^{-s}
\times
}
\\
&&
||(\sum_{i \geq (-\log_{2}l(P))\vee 0}(2^{is'}
{\rm osc}^{k-1}_p
f(x, i))^q)^{1/q}||_{L^p(P)} 
\\
&\leq&
 C
\sup_{ \mathcal{D} \ni Q \ni x_{0} }l(Q)^{-\sigma}
\sup_{ \mathcal{D} \ni P \subset 3Q\ }l(P)^{-s}\times
\\
&&
(||f||_{L^p(P)}+||(\sum_{i \geq (-\log_{2}l(P))\vee 0}
(2^{is'}d^{k}_if)^q)^{1/q})||_{L^p(P)}). 
\end{eqnarray*} 
}

\bigskip

\noindent
{\it Proof}.\ \ We will see that for $k|u| \leq 2^{-i}$,
\bigskip

$
|\triangle^{k}_{u}f(x)|\leq C( \sum_{e=0}^{k}|f(x+eu)-P_{B_{i}(x+eu)}f(x+eu)|)$

$
\leq C\sum_{e=0}^{k}\sum_{l\geq i} \Omega^{k-1}_pf(x+eu, l).
$
\bigskip

We consider a sequence for $i< \cdots < m \rightarrow \infty$,

$$
 B_{i}(x+eu) \supset \cdots \supset B_{m}(x+eu) \supset 
\cdots \rightarrow x+eu.
$$
Then we have 
\begin{eqnarray*}
\lefteqn
{
\frac{1}{|B_{m}|}\int_{B_{m}}|f-P_{B_{i}}f|\ dy \leq 
\frac{1}{|B_{m}|}\int_{B_{m}}|f-P_{B_{m}}f|\ dy
}
\\
&+& 
\frac{1}{|B_{m}|}\sum_{l= i+1}^m\int_{B_{m}}|P_{B_l}f-P_{B_{l-1}}f|\ dy
\\
&\leq&
 \frac{1}{|B_{m}|}\int_{B_{m}}|f-P_{B_{m}}f|\ dy+
C\sum_{l= i+1}^m\frac{1}{|B_{l}|}\int_{B_{l}}|f-P_{B_{l-1}}f|\ dy
\\
&\leq&
\frac{1}{|B_{m}|}\int_{B_{m}}|f-P_{B_{m}}f|\ dy
+
C\sum_{l= i}^m\frac{1}{|B_{l}|}\int_{B_{l}}|f-P_{B_{l}}f|\ dy. 
\end{eqnarray*}
Hence, we have 
\begin{eqnarray*}
\lefteqn{
|f(x+eu)-P_{B_{i}}(x+eu)| = 
\lim_{m\rightarrow \infty}\frac{1}{|B_{m}|}\int_{B_{m}}|f-P_{B_{i}}f|\ dy
}
\\
&\leq&
\lim_{m\rightarrow \infty}\frac{1}{|B_{m}|}\int_{B_{m}}|f-P_{B_{m}}f|\ dy
+C\sum_{l= i}^{\infty}\frac{1}{|B_{l}|}\int_{B_{l}}|f-P_{B_{l}}f|\ dy
\\
&\leq&
 C\sum_{l=i}^{\infty}\Omega_{p}^{k-1}f(x+eu,l).
\end{eqnarray*}
Therefore, we have  for a dyadic cube $P$ with $l(P)=2^{-j}$,
\begin{eqnarray*}
\lefteqn{
(\sum_{i \geq j\vee 0}(2^{is'}\sup_{k|u| \leq 2^{-i}}||\triangle_{u}^{k}f||
_{L^p(P)})^q)^{1/q}
}
\\
&\leq&
C(\sum_{i \geq j\vee 0}(2^{is'}\sum_{l \geq i}||\Omega_{p}^{k-1}f(x, l)||
_{L^p(3P)})^q)^{1/q}
\\
&\leq&
C(\sum_{i \geq j\vee 0}(2^{is'}||\Omega_{p}^{k-1}f(x,i)||_{L^p(3P)})^q)^{1/q}
\\
&\leq&
 C(\sum_{i \geq j\vee 0}(2^{is'}||{\rm osc}_{p}^{k-1}f(x,i)||_{L^p(3P)})^q)^{1/q}
\\
&\leq&
C(|x_{0}-x_{P}|+2^{-j})^{\sigma}2^{-js}\sup_{ \mathcal{D} \ni Q \ni x_{0} }l(Q)^{-\sigma}\sup_{ \mathcal{D} \ni P \subset 3Q\ }l(P)^{-s}\times
\\
&&
(\sum_{i \geq j\vee 0}(2^{is'}||{\rm osc}_{p}^{k-1}f(x,i)||_{L^p(P)})^q)^{1/q}
\end{eqnarray*}
by using Hardy's inequality if $s' > 0$. This completes the proof of the first half of (i).

Next, we will  prove the last half of (i).

We consider a function $\theta \in \mathcal{S}$ such that supp\ $\theta \subset \{k|u| \leq 1 \}$ and $\int \theta(u) \ du =1$. We put 
$$
h_{i}(x) = \int (f(x)- \triangle^{k}_u f(x))\theta_{i}(u) \ du
$$
where $\theta_{i}(u)=2^{ni}\theta(2^iu)$. 
We claim that
\bigskip

$
{\rm osc}^{k-1}_pf(x,i) \leq C (\frac{1}{|B_{i}(x)|}\int_{B_{i}(x)}
|d^{k}_{i}f(y)|^p\ dy)^{1/p}
+
C{\rm osc}^{k-1}_{p}h_{i}(x,\ i).
$

\bigskip

We see that
\begin{eqnarray*}
\lefteqn{
{\rm osc}^{k-1}_pf(x,i)\sim \Omega^{k-1}_pf(x,i)=
}
\\
&&
(\frac{1}{|B_{i}(x)|}\int_{B_{i}(x)}|f(y)-P_{B_{i}(x)}f(y)|^p\ dy)^{1/p}
\\
&\leq&
(\frac{1}{|B_{i}(x)|}\int_{B_{i}(x)}|f(y)-h_{i}(y)|^p\ dy)^{1/p}
\\
&+&
(\frac{1}{|B_{i}(x)|}\int_{B_{i}(x)}|h_{i}(y)-P_{B_{i}(x)}h_{i}(y)|^p\ dy)^{1/p}
\\
&+&
(\frac{1}{|B_{i}(x)|}\int_{B_{i}(x)}|P_{B_{i}(x)}h_{i}(y)-P_{B_{i}(x)}f(y)|^p\ dy)^{1/p}
\\
&\leq& 
C(\frac{1}{|B_{i}(x)|}\int_{B_{i}(x)}|f(y)-h_{i}(y)|^p\ dy)^{1/p}
\\
&+&
C(\frac{1}{|B_{i}(x)|}\int_{B_{i}(x)}|h_{i}(y)-P_{B_{i}(x)}h_{i}(y)|^p\ dy)^{1/p}
\\
&\leq&
 C(\frac{1}{|B_{i}(x)|}\int_{B_{i}(x)}(\int_{k|u| \leq 2^{-i}}
|\triangle^{k}_uf(y)||\theta_{i}(u)|\ du)^p\ dy)^{1/p}
\\
&&
+ C\Omega^{k-1}_ph_{i}(x,i)
\\
&\leq&
C(\frac{1}{|B_{i}(x)|}\int_{B_{i}(x)}
|d^{k}_{i}f(y)|^p\ dy)^{1/p}+
C{\rm osc}^{k-1}_ph_{i}(x,i).
\end{eqnarray*}
Next, we will estimate  ${\rm osc}^{k-1}_ph_{i}(x,i)$. We consider the $(k-1)$th Taylor polynomial $q(x)$ of $h_{i}$ at $x$. Then we have 
\begin{eqnarray*}
\lefteqn{h_{i}(y)-q(y)
}
\\
&=&
\int_{0}^{1}\sum_{|\beta|=k}\frac{k}{\beta!}\partial^{\beta}h_{i}(x+t(y-x))
(x-y)^{\beta}(1-t)^{k-1}\ dt
\\
&=&
\int_{0}^{1}\sum_{|\beta|=k}\frac{k}{\beta!}\int \sum^{k}_{m=1}
\left(
\begin{array}{c}
k \\
m
\end{array}
\right)
(-1)^{k-m}\partial^{\beta}f(x+t(y-x)+mu)\times
\\
&&
\theta_{l}(u)\ du (x-y)^{\beta}(1-t)^{k-1}\ dt
\\
&=&
\int_{0}^{1}\sum_{|\beta|=k}\frac{k}{\beta!}\int \sum^{k}_{m=1}
\left(
\begin{array}{c}
k \\
m
\end{array}
\right)
(-1)^{k-m}m^{k}f(x+t(y-x)+m2^{-i}u)\times
\\
&&
\partial^{\beta}\theta(u)\ du (x-y)^{\beta}(1-t)^{k-1}\ dt.
\end{eqnarray*}
Hence, we see by using Minkowski's inequality
\begin{eqnarray*}
\lefteqn{||{\rm osc}_p^{k-1}h_{i}(x, i)||_{L^p(P)}\leq ||(\frac{1}{|B_{i}(x)|}\int_{B_{i}(x)}
|h_{i}(y)-q(y)|^p\ dy)^{1/p}||_{L^p(P)}
}
\\
&\leq&
C(\int_{P}\frac{1}{|B_{i}(x)|}\int_{B_{i}(x)}(\int^1_0\int_{k|u| \leq 1}
\sum_{m=1}^{k}|f(x+t(y-x)+m2^{-i}u)|\times
\\
&&
|\partial^{\beta}\theta(u)| du  |x-y|^{k}(1-t)^{k-1}\ dt)^p\ dydx)^{1/p}
\\
\lefteqn{
\leq
C\int^1_0\int_{k|u| \leq 1}\sum_{m=1}^{k}(\frac{1}{|B_{i}(x)|}\int_{B_{i}(0)}
\int_{P}|f(x+ty+m2^{-i}u)|^p\ dxdy)^{1/p} \times 
}
\\
&&
2^{-ik}(1-t)^{k-1}\ dudt
\\
\lefteqn{
\leq
C\int^1_0\int_{k|u| \leq 1}\sum_{m=1}^{k}(\frac{1}{|B_{i}(x)|}\int_{B_{i}(0)}
\int_{P+ty+m2^{-i}u}|f(x)|^p\ dxdy)^{1/p} \times 
}
\\
&&
2^{-ik}(1-t)^{k-1}\ dudt
\\
&\leq&
C2^{-ik}(\int_{5P}|f(x)|^p\ dx)^{1/p} \leq C2^{-ik}||f||_{L^p(5P)}.
\end{eqnarray*}
Moreover, we have 
\begin{eqnarray*}
\lefteqn{
||(\frac{1}{|B_{i}(x)|}\int_{B_{i}(x)}
|d^{k}_{i}f(y)|^p\ dy)^{1/p}||_{L^p(P)}
}
\\
&\leq&
C
(\int_{P}(\frac{1}{|B_{i}(x)|}\int_{B_{i}(0)}
|d^{k}_{i}f(x+y)|^p\ dydx)^{1/p}
\\
&\leq&
 C(\frac{1}{|B_{i}(x)|}\int_{B_{i}(0)}\int_{P+y}
|d^{k}_{i}f(x)|^p\ dxdy)^{1/p}
\\
&\leq&
C(\int_{3P}|d^{k}_{i}f(x)|^p\ dx)^{1/p}\leq C||d^{k}_{i}f||_{L^p(3P)}.
\end{eqnarray*}
Thus, we have for a dyadic cube $P$ with $l(P)=2^{-j}$
\begin{eqnarray*}
\lefteqn{
(\sum_{i \geq j\vee 0}(2^{is'}||{\rm osc}^{k-1}_pf(x,i)||_{L^p(P)})^q)^{1/q}
}
\\
\lefteqn{
\leq
C(\sum_{i \geq j\vee 0}(2^{is'}||d^k_{i}f||_{L^p(3P)})^q)^{1/q}
+ 
C(\sum_{i \geq j\vee 0}2^{-i(k-s')q})^{1/q}||f||_{L^p(5P)}
}
\\
&\leq&
C(\sum_{i \geq j\vee 0}(2^{is'}||d^k_{i}f||_{L^p(3P)})^q)^{1/q}
+ 
C||f||_{L^p(5P)}
\\
\lefteqn{
\leq
 C(|x_{0}-x_{P}|+2^{-j})^{\sigma}2^{-js} \sup_{\mathcal{D} \ni Q \ni x_{0}}
l(Q)^{-\sigma}\sup_{ \mathcal{D} \ni P \subset 3Q}l(P)^{-s}\times
}
\\
&&
(\sum_{i \geq j\vee 0}(2^{is'}||d^{k}_{i}f||_{L^p(P)})^q)^{1/q}
\\
&+&
 C(|x_{0}-x_{P}|+2^{-j})^{\sigma}2^{-js} \sup_{\mathcal{D} \ni Q \ni x_{0}}
l(Q)^{-\sigma}\sup_{ \mathcal{D} \ni P \subset 3Q}l(P)^{-s}
||f||_{L^p(P)}
\end{eqnarray*}
if $k > s'$. The proof of (i) is complete. 
In the same way we can prove (ii).
\qed

\bigskip

\noindent
{\bf Theorem 6}. \ \ {\it 
{\rm (i)}\ \ Let $s',\ s,\ \sigma \in {\Bbb R}$ with $0 <s',\ \ 0\leq \sigma$, and let  $x_{0} \in \mathbb{R}^n$,
$ 1 \leq p \leq \infty,\ 0< q \leq \infty$. Let $k \in \mathbb{N}$ with $k > s'> 0$. Then we have following equivalences for $f \in \mathcal{S}'$ 
\begin{eqnarray*}
\lefteqn{
||f||_{A^s(B^{s'}_{pq})^{\sigma}_{x_{0}}}+\sup_{x_{0} \in Q}l(Q)^{-\sigma}\sup_{P \subset 3Q}l(P)^{-s}||f||_{L^p(P)}
}
\\
&\sim&
\sup_{x_{0} \in Q}l(Q)^{-\sigma}\sup_{P \subset 3Q}l(P)^{-s}
(||f||_{L^p(P)}
\\
&+&
(\sum_{i \geq (-\log_{2}l(P))\vee 0}(2^{is'}
\sup_{k|u|\leq 2^{-i}}||\triangle^{k}_{u}f||_{L^p(P)})^q)^{1/q})
\\
&\sim&
\sup_{x_{0}\in Q}l(Q)^{-\sigma}\sup_{P \subset 3Q}l(P)^{-s}
(||f||_{L^p(P)}
\\
&+&
(\sum_{i \geq (-\log_{2}l(P))\vee 0}(2^{is'}||{\rm osc}^{k-1}_{p}f||_{L^p(P)})^q)^{1/q})
\\
&\sim&
\sup_{x_{0}\in Q}l(Q)^{-\sigma}\sup_{P\subset 3Q}l(P)^{-s}
(||f||_{L^p(P)}
\\
&+&
(\sum_{i \geq (-\log_{2}l(P))\vee 0}(2^{is'}||d^{k}_{i}f||_{L^p(P)})^q)^{1/q}).
\end{eqnarray*}

{\rm (ii)}\ \ 
Let $s,\ s',\ \sigma \in {\Bbb R}$ with $0 <s',\ \ 0\leq \sigma,\ $  $x_{0} \in \mathbb{R}^n$,
$ 1 \leq p < \infty,\ 1\leq q \leq \infty$. Let  $k \in \mathbb{N}$ with $k > s'> 0$. Then we have following equivalences for $f \in \mathcal{S}'$ 

\begin{eqnarray*}
\lefteqn{
||f||_{A^s(F^{s'}_{pq})^{\sigma}_{x_{0}}}+\sup_{x_{0} \in Q}l(Q)^{-\sigma}\sup_{P \subset 3Q}l(P)^{-s}||f||_{L^p(P)}
}
\\
&\sim&
\sup_{x_{0} \in Q}l(Q)^{-\sigma}\sup_{P \subset 3Q}l(P)^{-s}
(||f||_{L^p(P)}
\\
&+&
||(\sum_{i \geq (-\log_{2}l(P))\vee 0}(2^{is'}
\sup_{k|u|\leq 2^{-i}}|\triangle^{k}_{u}f|)^q)^{1/q}||_{L^p(P)})
\\
&\sim&
\sup_{x_{0}\in Q}l(Q)^{-\sigma}\sup_{P \subset 3Q}l(P)^{-s}
(||f||_{L^p(P)}
\\
&+&
||(\sum_{i \geq (-\log_{2}l(P))\vee 0}(2^{is'}{\rm osc}^{k-1}_{p}f)^q)^{1/q}||_{L^p(P)})
\\
&\sim&
\sup_{x_{0}\in Q}l(Q)^{-\sigma}\sup_{P\subset 3Q}l(P)^{-s}
(||f||_{L^p(P)}
\\
&+&
||(\sum_{i \geq (-\log_{2}l(P))\vee 0}(2^{is'}d^{k}_{i}f)^q)^{1/q}||_{L^p(P)}).
\end{eqnarray*}
}

\bigskip

\noindent
{\it Proof}.\ \
(i)  It suffices to prove the first part of (i) by 
 Lemma 5. We consider the Littlewood-Paley decomposition 
$f=S_{i}f+ \sum_{l > i}f*\phi_{l}$. 
Then we have 
for $k|u| \leq 2^{-i}$ and a dyadic cube $P$ with $l(P)=2^{-j}$ , $i \geq j$
\begin{eqnarray*}
\lefteqn{
||\triangle^k_uf||_{L^p(P)}
\leq
||\triangle^k_u(f-S_{i}f)||_{L^p(P)}+||\triangle^k_uS_{i}f||_{L^p(P)}
}
\\
&\leq&
C\sum_{l > i}||\triangle^k_u(f*\phi_{l})||_{L^p(P)}+C||\triangle^k_uS_{i}f||_{L^p(P)}.
\end{eqnarray*}
We will  estimate $||\triangle^k_uS_{i}f||_{L^p(P)}$. 
Note the following formula 
$$
\triangle^k_uS_{i}f(x)=\int_{-\infty}^{\infty}\sum_{|\nu| =k}\frac{k!}{\nu!}
u^{\nu}\partial^{\nu}S_{i}f(x+\xi u)N_{k}(\xi)\ \ d\xi
$$
where $N_{k}$ is the B-spline of order $k$ (e.g. See [27]). Therefor we have for $k|u| \leq 2^{-i}$ 

\bigskip

$
||\triangle^k_uS_{i}f||_{L^p(P)} \leq C\sum_{|\nu|=k}|u|^k
||\partial^{\nu}S_{i}f||_{L^p(2P)}.
$

\bigskip

Next, we will  estimate $||\partial^{\nu}S_{i}f||_{L^p(2P)}$: 
\begin{eqnarray*}
\lefteqn{
||\partial^{\nu}S_{i}f||_{L^p(2P)}=||\int f(x-2^{-i}y) \
\partial^{\nu}\phi_{0}(y)\ dy||_{L^p(2P)}
}
\\
&\leq&
C\int(\int_{2P+2^{-i}y}|f(x)|^p\ dx)^{1/p} |\partial^{\nu}\phi_{0}(y)|\ dy
\\
&\leq&
 C2^{-js}\int(|x_{0}-x_{P}|+2^{-j}(1+|y|))^{\sigma} |\partial^{\nu}\phi_{0}(y)|\ dy
\\
&\times&
\sup_{x_{0}\in Q}l(Q)^{-\sigma}\sup_{P \subset 3Q}l(P)^{-s}||f||_{L^p(P)}\\
&\leq&
 C(|x_{0}-x_{P}|+2^{-j})^{\sigma}2^{-js}
\sup_{x_{0}\in Q}l(Q)^{-\sigma}\sup_{P \subset 3Q}l(P)^{-s}||f||_{L^p(P)}.
\end{eqnarray*}
Hence, we have 
\begin{eqnarray*}
\lefteqn{
||\triangle^k_uf||_{L^p(P)}
}
\\
&&
\leq
C\sum_{l > i}||\triangle^k_u(f*\phi_{l})||_{L^p(P)}
\\
&&
+
C(|x_{0}-x_{P}|+2^{-j})^{\sigma}2^{-js}2^{-ik}\sup_{x_{0}\in Q}l(Q)^{-\sigma}\sup_{P \subset 3Q}l(P)^{-s}||f||_{L^p(P)}.
\end{eqnarray*}
Moreover, we obtain by using Hardy's inequality if $s' >0$
\begin{eqnarray*}
\lefteqn{
(\sum_{i \geq j\vee 0}(2^{is'}\sup_{k|u| \leq 2^{-i}}||\triangle^k_uf||_{L^p(P)})^q)^{1/q}
}
\\
&\leq&
C(\sum_{i \geq j\vee 0}(2^{is'}\sum_{l > i}||\triangle^k_u(f*\phi_{l})||_{L^p(P)})^q)^{1/q}
\\
&+&
C(\sum_{i \geq j\vee 0}(2^{-i(k-s')}(|x_{0}-x_{P}|+2^{-j})^{\sigma}2^{-js}
\\
&\times&
\sup_{x_{0}\in Q}l(Q)^{-\sigma}\sup_{P \subset 3Q}l(P)^{-s}||f||_{L^p(P)})^q)^{1/q}
\\
&\leq&
C(\sum_{i > j\vee 0}(2^{is'}||\triangle^k_u(f*\phi_{i})||_{L^p(P)})^q)^{1/q}
\\
&+&
C(|x_{0}-x_{P}|+2^{-j})^{\sigma}2^{-js}
\sup_{x_{0}\in Q}l(Q)^{-\sigma}\sup_{P \subset 3Q}l(P)^{-s}||f||_{L^p(P)}.
\end{eqnarray*}
This implies that 
\begin{eqnarray*}
\lefteqn{
\sup_{x_{0}\in Q}l(Q)^{-\sigma}\sup_{P \subset 3Q}l(P)^{-s}
(\sum_{i \geq j\vee 0}(2^{is'}\sup_{k|u| \leq 2^{-i}}||\triangle^k_uf||_{L^p(P)})^q)^{1/q}
}
\\
&\leq&
C\sup_{x_{0}\in Q}l(Q)^{-\sigma}\sup_{P \subset 3Q}l(P)^{-s}
(\sum_{i > j\vee 0}(2^{is'}||f*\phi_{i}||_{L^p(P)})^q)^{1/q}
\\
&+&
C\sup_{x_{0}\in Q}l(Q)^{-\sigma}\sup_{P \subset 3Q}l(P)^{-s}||f||_{L^p(P)}
\\
&\leq&
C||f||_{A^s(B^{s'}_{pq})^{\sigma}_{x_{0}}}
+C\sup_{x_{0}\in Q}l(Q)^{-\sigma}\sup_{P \subset 3Q}l(P)^{-s}||f||_{L^p(P)}.
\end{eqnarray*} 
We will show the converse statement. 
It is easy to see that there exist $\phi^m \in \mathcal{S}\ \ m=1, \cdots n$ such that 
$ \phi=\sum_{m=1}^n\triangle^k_{ce_{m}}\phi^m$ for enough small $c$ 
where $e_{1},\cdots,e_{n}$ are the canonical basis vectors in $\mathbb{R}^n$. 
Then we have for $i \in \mathbb{N}$
$$
f*\phi_{i}=\sum_{m=1}^nf*\triangle^k_{c2^{-i}e_{m}}\phi^m_{i}= 
\sum_{m=1}^n\triangle^k_{c2^{-i}e_{m}}f*\phi^m_{i}. 
$$
Therefore, we have for a dyadic cube $P$ with $l(P)=2^{-j}$ and $i \geq j$
\begin{eqnarray*}
\lefteqn{
||f*\phi_{i}||_{L^p(P)} 
}
\\
&\leq&
C||\sum_{m=1}^n\triangle^k_{c2^{-i}e_{m}}f*\phi^m_{i}||_{L^p(P)}
\\
&\leq&
C\int \sum_{m=1}^n(\int_{P+2^{-i}y}|\triangle^k_{c2^{-i}e_{m}}f(x)|^p\ dx)^{1/p}
|\phi^m(y)|\ dy
\\
&\leq&
C\int\sum_{m=1}^n (\int_{P+2^{-i}y}\sup_{k|u| \leq 2^{-i}}|\triangle_{u}^kf(x)|^p\ dx)^{1/p}|\phi^m(y)|\ dy.
\end{eqnarray*} 
Hence, we have if $l(P) < 1$
\begin{eqnarray*}
\lefteqn{
(\sum_{i \geq j}(2^{is'}||f*\phi_{i}||_{L^p(P)})^q)^{1/q}
}
\\
&\leq&
C(|x_{0}-x_{P}|+2^{-j})^{\sigma}2^{-js}\times
\\
&&
\sup_{x_{0}\in Q}l(Q)^{-\sigma}\sup_{P \subset 3Q}l(P)^{-s}
(\sum_{i \geq j}(2^{is'}\sup_{k|u| \leq 2^{-i}}||\triangle^k_u f||_{L^p(P)})^q)^{1/q}
\end{eqnarray*} 
and if $l(P) \geq 1$
\begin{eqnarray*}
\lefteqn{
(\sum_{i \geq 0}(2^{is'}||f*\phi_{i}||_{L^p(P)})^q)^{1/q}
}
\\
&\leq&
(\sum_{i > 0}(2^{is'}||f*\phi_{i}||_{L^p(P)})^q)^{1/q}
+||f*\phi_{0}||_{L^p(P)}
\\
&\leq&
C(|x_{0}-x_{P}|+2^{-j})^{\sigma}2^{-js}
\sup_{x_{0}\in Q}l(Q)^{-\sigma}\sup_{P \subset 3Q}l(P)^{-s}\times
\\
&&
(\sum_{i >0}(2^{is'}\sup_{k|u| \leq 2^{-i}}||\triangle^k_u f||_{L^p(P)})^q)^{1/q}
\\
&+&
C(|x_{0}-x_{P}|+2^{-j})^{\sigma}2^{-js}
\sup_{x_{0}\in Q}l(Q)^{-\sigma}\sup_{P \subset 3Q}l(P)^{-s}
|| f||_{L^p(P)}.
\end{eqnarray*} 
Thus, we have 
\begin{eqnarray*}
\lefteqn{
\sup_{x_{0}\in Q}l(Q)^{-\sigma}\sup_{P \subset 3Q}l(P)^{-s}(\sum_{i \geq 0}(2^{is'}||f*\phi_{i}||_{L^p(P)})^q)^{1/q}
}
\\
&\leq&
C\sup_{x_{0}\in Q}l(Q)^{-\sigma}\sup_{P \subset 3Q}l(P)^{-s}\times
\\
&&
(\sum_{i \geq 0}(2^{is'}\sup_{k|u| \leq 2^{-i}}||\triangle^k_u f||_{L^p(P)})^q)^{1/q}
\\
&+&
\sup_{x_{0}\in Q}l(Q)^{-\sigma}\sup_{P \subset 3Q}l(P)^{-s}
|| f||_{L^p(P)}.
\end{eqnarray*} 
This completes the proof of Theorem 6 (i). 
In the same way we can prove (ii).

\qed

\bigskip

\noindent

\section*{\bf\large Acknowlegments} 
The author would like to thank Prof. Yoshihiro Sawano for 
his encouragement and  many helpful remarks. Furthermore, the author would like to thank  the referee for most helpful advices and corrections.


\vskip 1cm \footnotesize
\begin{flushleft}
Koichi~Saka \\
Department of Mathematics\\
Akita University, \\
010-8502 Akita, Japan\\

E-mail:  sakakoichi@gmail.com
\end{flushleft}

 \end{document}